\magnification=\magstep1
\input amstex
\documentstyle{amsppt}
\pagewidth{30pc}\pageheight{48pc}
\input xy
\xyoption{all}

\loadeusm
\define\Z{\Bbb Z}

\topmatter
\title Stably dualizable groups \endtitle
\author John Rognes \endauthor
\date September 28th 2005 \enddate
\address Department of Mathematics, University of Oslo, Norway \endaddress
\email rognes\@math.uio.no \endemail
\abstract
We extend the duality theory for topological groups from the classical
theory for compact Lie groups, via the topological study by J.~R.~Klein
\cite{Kl01} and the $p$-complete study for $p$-compact groups by
T.~Bauer \cite{Ba04}, to a general duality theory for stably dualizable
groups in the $E$-local stable homotopy category, for any spectrum $E$.
The principal new examples occur in the $K(n)$-local category, where
the Eilenberg--Mac\,Lane spaces $G = K(\Z/p, q)$ are stably dualizable
and nontrivial for $0 \le q \le n$.

We show how to associate to each $E$-locally stably dualizable group~$G$
a stably defined representation sphere $S^{adG}$, called the dualizing
spectrum, which is dualizable and invertible in the $E$-local category.
Each stably dualizable group is Atiyah--Poincar{\'e} self-dual
in the $E$-local category, up to a shift by $S^{adG}$.  There are
dimension-shifting norm- and transfer maps for spectra with $G$-action,
again with a shift given by $S^{adG}$.  The stably dualizable group $G$
also admits a kind of framed bordism class $[G] \in \pi_*(L_E S)$, in
degree $\dim_E(G) = [S^{adG}]$ of the $Pic_E$-graded homotopy groups of
the $E$-localized sphere spectrum.
\endabstract

\toc
\head 1. Introduction \endhead
\subhead 1.1. The symmetry groups of stable homotopy theory \endsubhead
\subhead 1.2. Algebraic localizations and completions \endsubhead
\subhead 1.3. Chromatic localizations and completions \endsubhead
\subhead 1.4. Applications \endsubhead
\head 2. The dualizing spectrum \endhead
\subhead 2.1. The $E$-local stable category \endsubhead
\subhead 2.2. Dualizable spectra \endsubhead
\subhead 2.3. Stably dualizable groups \endsubhead
\subhead 2.4. $E$-compact groups \endsubhead
\subhead 2.5. The dualizing and inverse dualizing spectra \endsubhead
\head 3. Duality theory \endhead
\subhead 3.1. Poincar{\'e} duality \endsubhead
\subhead 3.2. Inverse Poincar{\'e} duality \endsubhead
\subhead 3.3. The Picard group \endsubhead
\head 4. Computations \endhead
\subhead 4.1. A spectral sequence for $E$-homology \endsubhead
\subhead 4.2. Morava $K$-theories \endsubhead
\subhead 4.3. Eilenberg--Mac\,Lane spaces \endsubhead
\head 5. Norm and transfer maps \endhead
\subhead 5.1. Thom spectra \endsubhead
\subhead 5.2. The norm map and Tate cohomology \endsubhead
\subhead 5.3. The $G$-transfer map \endsubhead
\subhead 5.4. $E$-local homotopy classes \endsubhead
\head {\ } References \endhead
\endtoc
\endtopmatter

\define\F{\Bbb F}
\define\G{\Bbb G}
\define\Hom{\operatorname{Hom}}
\define\K{\eusm K}

\define\M{\eusm M}
\define\Pic{\operatorname{Pic}}

\define\Q{\Bbb Q}
\define\Tor{\operatorname{Tor}}
\define\Tot{\operatorname{Tot}}

\define\trf{\operatorname{trf}}
\redefine\D{\eusm D}
\redefine\L{\eusm L}
\redefine\phi{\varphi}

\document

\head 1. Introduction \endhead

\subhead 1.1. The symmetry groups of stable homotopy theory \endsubhead

Compact Lie groups occur naturally as the symmetry groups of geometric
objects, e.g.~as the isometry groups of Riemannian manifolds \cite{MS39}.
Such geometric objects can usefully be viewed as equivariant objects,
i.e., as a spaces with an action by a Lie group.  The homotopy theory of
such equivariant spaces is quite well approximated by the corresponding
stable equivariant homotopy theory, which in its strong ``genuine''
form relies, already in its construction, on the good representation
theory for actions by Lie groups on finite-dimensional vector spaces.

As a first example of a useful stable result, consider the Adams
equivalence $Y/G \simeq (\Sigma^{-adG} Y)^G$ of \cite{LMS86, II.7}.  Here
$Y$ is any free $G$-spectrum, $adG$ denotes the adjoint representation
of $G$ on its Lie algebra and $\Sigma^{-adG} Y$ is the stably defined
desuspension of $Y$ with respect to this $G$-representation.

As a second example, Atiyah duality \cite{At61} asserts that if $M$
is a smooth closed manifold with stable normal bundle $\nu$, the
functional (Spanier--Whitehead) dual $DM_+ = F(M_+, S)$ of $M_+$ is
equivalent to the Thom spectrum $Th(\nu \downarrow M)$.  When $M = G$
is a compact Lie group, and thus parallelizable, we can write this as
a stable Poincar{\'e} duality equivalence $DG_+ \simeq Th(\epsilon^{-n}
\downarrow G) = \Sigma^{-n} \Sigma^\infty(G_+)$.  But $G$ acts on itself
both from the left and the right, and the bi-equivariant form of this
equivalence takes the more precise form
$$
DG_+ \wedge S^{adG} \simeq \Sigma^\infty(G_+)
$$
where $G$ acts by conjugation from the left on the one-point
compactification $S^{adG}$ of the adjoint representation and trivially
from the right.  See Theorem~3.1.4 below.

As a third example, the left-invariant framing of an $n$-dimensional
compact Lie group $G$ gives it an associated stably framed cobordism
class $[G]$ in $\Omega^{fr}_n \cong \pi_n(S)$, the $n$-th stable stem.
For example $[S^1] = \eta \in \pi_1(S)$ realizes the stable class of the
Hopf fibration $\eta \: S^3 \to S^2$.  It is of interest to see which
stable homotopy classes actually occur in this way \cite{Os82}.

The formulation of these three results may appear to require that $G$
admits a geometric representation theory, with tangent spaces, adjoint
representations, etc., but in fact much less is required, and that is
the main thrust of the present article.

\subhead 1.2. Algebraic localizations and completions \endsubhead

Homotopy-theoretically, the main properties of compact Lie groups are
(i) that they are compact manifolds, hence admit the structure of a
finite CW complex, and (ii) that they are topological groups, hence are
homotopy equivalent to loop spaces.  Browder \cite{Brd61, 7.9} showed
that all finite H-spaces are Poincar{\'e} complexes, and recently Bauer,
Kitchloo, Notbohm and Pedersen \cite{BKNP04} showed that all finite loop
spaces indeed are homotopy equivalent to manifolds.  However, there are
examples of finite loop spaces that are not even rationally equivalent
to Lie groups \cite{ABGP04}.

A standard method in homotopy theory, and a key ingredient in
\cite{BKNP04}, is the possibility to study homotopy types locally, say
by Bousfield localization with respect to a homology theory \cite{Bo75},
\cite{Bo79}, or completion in the sense of Bousfield and Kan \cite{BK72}.
For later generalization, we recall that the Bousfield--Kan $p$-completion
functor is derived from a monad $(\F_p(-), \mu, \eta)$, where $X \mapsto
\F_p(X)$ is a particular endofunctor on spaces with $\pi_* \F_p(X) =
\widetilde H_*(X; \F_p)$, equipped with suitable associative and unital
natural transformations $\mu_X \: (\F_p \circ \F_p)(X) \to \F_p(X)$
and $\eta_X \: X \to \F_p(X)$.  For each space $X$ the monad generates
a cosimplicial space
$$
[q] \longmapsto (\F_p \circ \dots \circ \F_p)(X)
$$
(the $(q+1)$-fold iterate of the endofunctor), whose totalization defines
the Bousfield--Kan $p$-completion of $X$, see \cite{BK72, \S I.4},
which we denote by $X^\wedge_p$.  We say that $X$ is $p$-complete when
the natural map $X \to X^\wedge_p$ is a weak equivalence.  Each space
$\F_p(X)$ is a product of mod~$p$ Eilenberg--Mac\,Lane spaces, i.e.,
the spaces in the $\Omega$-spectrum of the mod~$p$ Eilenberg--Mac\,Lane
spectrum $H\F_p$, and $X^\wedge_p$ is constructed as a limit of such
spaces.

In the Bousfield--Kan $p$-complete category, the local incarnations
of finite loop spaces are the {\it $p$-compact groups\/} of Dwyer and
Wilkerson \cite{DW94}.  These are the topological groups $G$ with finite
mod~$p$ homology $H_*(G; \F_p)$ and $p$-adically complete classifying
space $BG \simeq BG^\wedge_p$.  Dwyer and Wilkerson consider loop
spaces instead of topological groups, but any loop space is equivalent
to a topological group, and conversely, so there is no real distinction.
We shall prefer to work with topological groups $G$ and their classifying
spaces $BG$, rather than with loop spaces $\Omega Y$ and their deloopings
$Y$, since we shall make extensive use of group actions by $G$, which
would involve operad actions or coherence theory to formulate for
loop spaces.

We think of a connected compact Lie group $G$ as a geometric, integrally
defined object, which can be analyzed one rational prime $p$ at a time by
way of its homotopy-theoretic, locally defined $p$-compact pieces, namely
the $p$-compact groups $\Omega (BG)^\wedge_p$ obtained by $p$-completing
the classifying space $BG$ at $p$ and looping.  In this case, the
fact that the loop space of $(BG)^\wedge_p$ still has finite mod~$p$
homology follows from the convergence of the mod~$p$ Eilenberg--Moore
spectral sequence \cite{Dw74}.  In addition to these geometric examples,
there are also other ``exotic'' $p$-compact groups that only exist
locally, without the global origin of a compact Lie group, such as
the Dwyer--Wilkerson $2$-compact group $DI(4)$, for which $H^*(BDI(4);
\F_2)$ realizes the ring of rank~$4$ Dickson invariants \cite{DW93}.
There is no compact Lie group with this cohomology algebra.

In his Ph.D.~thesis, T.~Bauer \cite{Ba04} showed that for each $p$-compact
group~$G$ one can produce a $p$-complete stable replacement for the
adjoint representation sphere $S^{adG}$, for the purposes of $p$-complete
stable homotopy theory.  It suffices to work $G$-equivariantly in the
``naive'' sense, where the objects are spectra equipped with a $G$-action,
and the (weak) equivalences are $G$-equivariant maps that are stable
equivalences in the underlying non-equivariant category.  Bauer showed
that for a $p$-compact group $G$, analogous results to the Adams
equivalence and the Atiyah--Poincar{\'e} duality equivalences above hold,
with $S^{adG}$ reinterpreted as the dualizing spectrum $(\Sigma^{\infty}
G_+)^{hG} = F(EG_+, \Sigma^{\infty} G_+)^G$ of W.~Dwyer (unpublished)
and J.~R.~Klein \cite{Kl01}, but formed in the $p$-complete category.
Bauer also showed that a $p$-compact group $G$ has the analogue of
a framed bordism class $[G]$ in $\pi_*(S^\wedge_p)$.  For example,
the Sullivan spheres (see Example~2.4.4) are examples of $p$-compact
groups for $p$ odd, and their framed bordism classes represent the
generators $\alpha_1 \in \pi_{2p-3}(S^\wedge_p)$.  We shall generalize
the constructions and results of Bauer to other, topologically defined,
local categories.

\subhead 1.3. Chromatic localizations and completions \endsubhead

In stable homotopy theory it is well-known, following Ravenel \cite{Ra84},
that for each prime $p$ it is possible to interpolate in infinitely many
``chromatic'' stages between the algebraically $p$-local and rational
situations, through a tower of Bousfield localizations with respect to
the homology theories represented by the Johnson--Wilson spectra $E(n)$,
for $n\ge0$.  Furthermore, one can isolate the individual strata of
this filtration by way of the Bousfield localization with respect to the
Morava $K$-theory spectra $K(n)$, for $n\ge1$.  In a precise sense, these
$K(n)$-local strata are the finest non-trivial localizations possible
\cite{HSt99, 7.5}.  We therefore have the possibility of analyzing the
stable images of compact Lie groups, or $p$-compact groups, in even finer
detail than that offered by the algebraic localizations, by working in
the $p$-primary $K(n)$-local category for one prime $p$ and one natural
number $n$, thereby focusing only on the ``$p$-primary $v_n$-periodic''
parts of their homotopy theory.

The topological groups $G$ that have the finiteness property that
$K(n)_*(G)$ is finite in each degree will be called {\it $K(n)$-locally
stably dualizable\/} groups.  Among these we can single out the {\it
$K(n)$-compact groups\/} as those whose classifying space $BG$ is a
$K(n)$-complete space, in the sense of Bendersky, Curtis and Miller
\cite{BCM78} and \cite{BT00, \S2}, which will be recalled in Section~2.4
below.  These are the precise analogues, to the eyes of Morava $K$-theory,
of the $p$-compact groups, to the eyes of mod~$p$ homology.  However, the
real work in this paper applies to all stably dualizable groups.  There
are examples of even more ``exotic'' $K(n)$-locally stably dualizable
groups than the $p$-compact ones, that only exist $K(n)$-locally for
some~$n$.  The simplest, abelian, instances of this are provided by the
Eilenberg--Mac\,Lane spaces $G = K(\pi, q)$, e.g.~for $\pi = \Z/p$, $0
\le q \le n$ \cite{RaW80, 9.2}, which for $q \ge 1$ have infinite mod~$p$
homology and thus are never $p$-compact.

In this paper we generalize the duality theory of Lie groups and
of $p$-compact groups by Klein and Bauer, to show that also for a
$K(n)$-locally stably dualizable group $G$, the {\it dualizing spectrum}
$$
S^{adG} = L_{K(n)} (\Sigma^\infty G_+)^{hG}
$$
formed in the $K(n)$-local stable category has the properties that
make it a stable substitute for the adjoint representation sphere of
a compact Lie group.  Here the $G$-homotopy fixed points are formed
with respect to the standard right $G$-action on $\Sigma^\infty G_+$.
The dualizing spectrum $S^{adG}$ is a {\it dualizable\/} and {\it
invertible\/} spectrum in the $K(n)$-local category, cf.~Theorem~3.3.4,
which means that it has an equivalence class
$$
[S^{adG}] \in \Pic_{K(n)}
$$
in the $K(n)$-local Picard group \cite{HMS94}.  In particular, suspending
(smashing) by $S^{adG}$ is an invertible self-equivalence of the
$K(n)$-local category.  The $K(n)$-local smash inverse of $S^{adG}$
is its functional dual $S^{-adG} = DS^{adG} = F(S^{adG}, L_{K(n)}S)$.
See Propositions~2.5.7 and~3.2.3.

We show that there is a natural {\it norm map}
$$
N \: (X \wedge S^{adG})_{hG} \to X^{hG}
$$
for any spectrum $X$ with $G$-action, which is a $K(n)$-local equivalence
under slightly different conditions on $X$ than those of the Adams
equivalence.  Up to rewriting, it is the canonical colim/lim exchange
map for the $G$-homotopy orbit and $G$-homotopy fixed point constructions
on $X \wedge \Sigma^\infty G_+$.  See Theorem~5.2.5.  We also show that
there is an (implicitly $K(n)$-local) natural {\it Atiyah--Poincar{\'e}
duality equivalence}
$$
DG_+ \wedge S^{adG} \simeq \Sigma^\infty G_+ \,,
$$
which is $G$-equivariant from both the left and the right.  See
Theorem~3.1.4.  Finally, we combine the norm map $N \: BG^{adG} =
(S^{adG})_{hG} \to S^{hG} = D(BG_+)$ for $X = S$ with a bottom cell
inclusion $i \: S^{adG} \to BG^{adG}$ and the projection $p \: S^{hG}
\to S$ to obtain a natural map
$$
pNi \: S^{adG} \to S \,,
$$
representing a homotopy class
$$
[G] \in \pi_*(L_{K(n)}S)
$$
in the {\it $\Pic_{K(n)}$-graded homotopy groups\/} of the $K(n)$-local
sphere spectrum.  See Definition~5.4.1.  We informally think of this as
the $K(n)$-locally {\it framed bordism class\/} of $G$.

The results discussed up to now hold in a uniform manner in the
$E$-local stable category, for each fixed spectrum~$E$ and suitably
defined $E$-locally stably dualizable groups.  The terminology is
chosen so that $G$ is $E$-locally stably dualizable precisely when its
suspension spectrum $\Sigma^\infty G_+$ is dualizable in the $E$-local
stable category.  See Definition~2.3.1.  This is the generality in which
the main body of the paper is written.

In Chapter~4 we develop calculational tools to study $E$-locally stably
dualizable groups, mostly particular to the most local case $E = K(n)$.
The group structure on~$G$ makes $H = K(n)_*(G)$ a graded Frobenius
algebra over $R = K(n)_*$ (Proposition 4.2.4), for the $R$-dual $H^* =
K(n)^*(G)$ is a free graded $H$-module of rank~$1$.  There is a strongly
convergent homological spectral sequence of Eilenberg--Moore type
$$
E^2_{s,t} = \Tor_{s,t}^{H}(R, H^*)
\Longrightarrow K(n)^{-(s+t)}(S^{adG})
$$
(Proposition~4.1.1).  It collapses at the $E^2$-term to the line $s=0$,
and its dual identifies $K(n)_*(S^{adG})$ with the $H^*$-comodule
primitives $P_{H^*}(H) \cong \Hom_H(H^*, R)$ in $\Hom_R(H^*, R) \cong H
= K(n)_*(G)$ (Theorem~4.2.6).  For example, when $G = K(\Z/p, n)$ is
viewed as a $K(n)$-locally stably dualizable group, it follows that
$[G] \: S^{adG} \to S$ is an equivalence in the $K(n)$-local category
(Example 5.4.6), so that the Atiyah--Poincar{\'e} duality equivalence
takes the untwisted form
$$
F(K(\Z/p, n)_+, L_{K(n)}S) \simeq L_{K(n)} \Sigma^{\infty} K(\Z/p, n)_+
\,.
$$

\subhead 1.4. Applications \endsubhead

It is conceivable that more invertible spectra in the $K(n)$-local
category can be constructed in the form $S^{adG}$ for $K(n)$-locally
stably dualizable groups $G$, than just the localized integer sphere
spectra $L_{K(n)} \Sigma^d S$ for $d \in \Z$.  There are no such
examples in the $p$-complete setting, but the $K(n)$-local Picard group
is more subtle.  Likewise, it is conceivable that the associated homotopy
classes $[G] \in \pi_*(L_{K(n)}S)$ can realize more homotopy classes than
those that appear from Lie groups and $p$-compact groups.  However, so
far we have mostly studied the abelian examples of $K(n)$-locally stably
dualizable groups given by Eilenberg--Mac\,Lane spaces, where this added
potential is not realized.  We think of these abelian groups as playing
the role analogous to tori in the theory of compact Lie groups, and hope
to develop a richer supply of non-abelian examples in joint work with
Bauer, cf.~Remark~2.4.8.

This work was simultaneously motivated by the author's formulation
\cite{Rog:g} of Galois theory of $E$-local commutative $S$-algebras.
If $A \to B$ is an $E$-local $G$-Galois extension there is a useful
norm equivalence $N \: (B \wedge S^{adG})_{hG} \to B^{hG}$, with $A
\simeq B^{hG}$.  For finite groups $G$ this follows as in \cite{Kl01},
but the natural generality for the theory appears to be to allow
topological Galois groups $G$ that are $E$-locally stably dualizable,
as considered here.  The constructions in Chapters~3 and~5 of the
present paper then find applications in the cited Galois theory.

\comment
Along the same lines, it is desirable to cast the theory of pro-Galois
extensions in a topological context.  For example, the height~$n$
Morava stabilizer group $\G_n$ is the pro-Galois group of the extension
$L_{K(n)}S \to E_n$, and it is a $p$-adic duality group in the sense that
the continuous group cohomology $H^*(\G_n; \Z_p[[\G_n]])$ is isomorphic to
$\Z_p$ concentrated in degree $* = n^2$ \cite{St00, Prop.~5}.  We expect
that it will follow, in a fashion similar to \cite{Kl01, 10.7}, that
$\G_n$ is $p$-completely stably dualizable, with dualizing spectrum
homotopy equivalent to the $p$-completion of $S^{-n^2}$.  Likewise in
the $K(n)$-local category.
\endcomment

\subhead Acknowledgments \endsubhead

The author wishes to thank Tilman Bauer for discussions starting with
\cite{Ba04} and leading to the present paper, and the referee for very
constructive comments.  Part of this work was done while the author
was a member of the Isaac Newton Institute for Mathematical Sciences,
Cambridge, in the fall of 2002, and he wishes to thank the INI for its
hospitality and support.

\head 2. The dualizing spectrum \endhead

\subhead 2.1. The $E$-local stable category \endsubhead

As our basic model for spectra we shall take the bicomplete, bitensored
closed symmetric monoidal category $\M_S$ of {\it $S$-modules\/} from
\cite{EKMM97}.  The symmetric monoidal pairing is the smash product $X
\wedge Y$, the unit object is the sphere spectrum~$S$, and the internal
function object is the mapping spectrum $F(X, Y)$.  We write $DX = F(X,
S)$ for the functional dual.  For a based topological space $T$ we write
$X \wedge T = X \wedge \Sigma^\infty T$ and $F(T, X) = F(\Sigma^\infty T,
X)$ for the resulting bitensors.

Let $E$ be any $S$-module.  It induces the (generalized, reduced) homology
theory $E_*$ that takes an $S$-module $X$ to the graded abelian group
$E_*(X) = \pi_*(E \wedge X)$.  A map $f \: X \to Y$ of $S$-modules is
said to be an {\it $E$-equivalence\/} if the induced homomorphism $f_*
\: E_*(X) \to E_*(Y)$ is an isomorphism, and an $S$-module $Z$ is {\it
$E$-local\/} if for each $E$-equivalence $f \: X \to Y$ the induced
homomorphism $f^\# \: [Y, Z]_* \to [X, Z]_*$ is an isomorphism.

Let $\M_{S,E}$ be the full subcategory of $\M_S$ of {\it $E$-local
$S$-modules}.  There is a Bousfield localization functor $L_E \: \M_S
\to \M_{S,E}$ \cite{Bo79}, \cite{EKMM97, Ch.~VIII} that comes equipped
with a natural $E$-equivalence $X \to L_E X$ for each $S$-module $X$
(with $L_E X$ $E$-local).  Let $\D_S = \bar h \M_S$ be the homotopy
category of $\M_S$, i.e., the {\it stable category}, and let $\D_{S,E} =
\bar h \M_{S,E}$ be the homotopy category of $\M_{S,E}$, i.e., the {\it
$E$-local stable category}.  It is a stable homotopy category in the sense
of \cite{HPS97, 1.2.2}.  The induced $E$-localization functor $L_E \:
\D_S \to \D_{S,E}$ is left adjoint to the forgetful functor $\D_{S,E}
\to \D_S$.

The $E$-local category $\M_{S,E}$ inherits the structure of a bicomplete,
bitensored closed symmetric monoidal category from $\M_S$ by applying
$L_E$ to each construction formed in $\M_S$.  The symmetric monoidal
pairing takes $X$ and~$Y$ to $L_E(X \wedge Y)$, and the unit object is the
$E$-local sphere spectrum $L_E S$.  The internal function object $F(X,
Y)$ is already $E$-local when $Y$ is $E$-local, hence does not change
when $E$-localized.  In a similar fashion the (limits and) colimits in
$\M_{S,E}$ are obtained from those formed in $\M_S$ by applying the
$E$-localization functor, and likewise for tensors (and cotensors).

\example{Example 2.1.1}
We may take $E = S$, in which case every spectrum is $S$-local, $\M_{S,S}
= \M_S$ and the $S$-local stable category is the whole stable category.
\endexample

\example{Example 2.1.2}
For a fixed rational prime $p$ and number $0 \le n < \infty$ we may take
$E = E(n)$, the $n$-th $p$-primary Johnson--Wilson spectrum, with
$$
E(n)_* = \Z_{(p)}[v_1, \dots, v_n, v_n^{-1}] \,.
$$
When $n=0$, $E(0) = H\Q$ is the rational Eilenberg--Mac\,Lane spectrum and
$E$-equivalence means rational equivalence.  In each case $L_n = L_{E(n)}$
is a smashing localization, $L_n S$ is a commutative $S$-algebra and the
$E(n)$-local category $\L_n = \M_{S,E(n)}$, as studied in \cite{HSt99},
is equivalent to the category $\M_{L_n S}$ of $L_n S$-modules.  In this
case the forgetful functor $\M_{S, E(n)} \to \M_S$ preserves the symmetric
monoidal pairing, but not the unit object.
\endexample

\example{Example 2.1.3}
For each prime $p$ and number $0 \le n \le \infty$ we may alternatively
take $E = K(n)$, the $n$-th $p$-primary Morava $K$-theory spectrum.
When $n=0$, $K(0) = E(0) = H\Q$, as discussed above.  When $0 < n <
\infty$,
$$
K(n)_* = \F_p[v_n, v_n^{-1}]
$$
is a graded field, and $\K_n = \D_{S, K(n)}$ is the $K(n)$-local stable
category, again studied in \cite{HSt99}.  When $n=\infty$, $K(\infty) =
H\F_p$ and $\M_{S, H\F_p}$ is the category of $H\F_p$-local $S$-modules.
For a connective spectrum $X$, $H\F_p$-localization amounts to algebraic
$p$-completion.  For $0 < n \le \infty$ the forgetful functor to $\M_S$
neither preserves the symmetric monoidal pairing nor the unit object.
\endexample

\definition{Convention 2.1.4}
Hereafter we shall work entirely within the $E$-local category $\M_{S,E}$.
We refer to the objects of $\M_{S,E}$ as {\it $E$-local $S$-modules}, or
simply as {\it spectra}.  For brevity we shall write $X \wedge Y$ for the
smash product $L_E(X \wedge Y)$, $S$ for the sphere spectrum $L_ES$ and
$F(X, Y)$ for the function spectrum $L_E F(X, Y)$ within this category.
The same applies to the functional dual $DX$, limits, colimits, tensors
and cotensors, all of which then take values in $\M_{S,E}$.
\enddefinition

\subhead 2.2. Dualizable spectra \endsubhead

Following Dold and Puppe \cite{DP80}, Lewis, May and Steinberger
\cite{LMS86, III.1} observe that in any closed symmetric monoidal category
there are natural canonical maps $\rho \: X \to DDX$, $\nu \: F(X, Y)
\wedge Z \to F(X, Y \wedge Z)$ and $\wedge \: F(X, Y) \wedge F(Z, W) \to
F(X \wedge Z, Y \wedge W)$.  We follow Hovey and Strickland \cite{HSt99,
1.5} and say that a spectrum $X$ is {\it $E$-locally dualizable\/}
if the canonical map
$$
\nu \: DX \wedge X \to F(X, X)
$$
(in the special case $X = Z$, $Y = S$) is an equivalence in $\M_{S,E}$.
When the spectrum $E$ is clear from the context, we simply say that $X$
is dualizable.  Lewis {\it et al.} then show \cite{LMS86, III.1.2, 1.3}:

\proclaim{Lemma 2.2.1}
\roster
\item
The canonical map $\rho \: X \to DDX$ is an
equivalence if $X$ is dualizable.
\item
The canonical map $\nu \: F(X, Y) \wedge Z \to F(X, Y \wedge Z)$ is an
equivalence if $X$ or $Z$ is dualizable.
\item
The smash product map $\wedge \: F(X, Y) \wedge F(Z, W) \to F(X \wedge Z,
Y \wedge W)$ is an equivalence if $X$ and $Z$ are dualizable,
or if $X$ is dualizable and $Y = S$.
\endroster
\endproclaim

It follows that the function spectrum $F(X, Y)$ and smash product $X
\wedge Y$ are dualizable when $X$ and $Y$ are dualizable.  In particular,
$DX$ is dualizable when $X$ is dualizable.

\example{Example 2.2.2}
For $E = S$, a spectrum $X$ is dualizable if and only if it is stably
equivalent to a finite CW spectrum, i.e., if and only if $X \simeq
\Sigma^\infty \Sigma^d K$ for some finite CW complex $K$ and integer
$d \in \Z$.  See e.g.~\cite{May96, XVI.7.4} for a proof, although the
non-equivariant result must be older.
\endexample

\example{Example 2.2.3}
For $E = K(n)$ with $0 \le n \le \infty$, Hovey and Strickland \cite{HSt99,
8.6} show that a $K(n)$-local $S$-module $X$ is dualizable if and only
if $K(n)_*(X)$ is a finitely generated $K(n)_*$-module.  Note that
this includes the cases $n=0$ with $K(0) = H\Q$ and $n = \infty$
with $K(\infty) = H\F_p$.  In each case $K(n)_*$ is a graded field,
so $K(n)_*(X)$ will automatically be free.
\endexample

\proclaim{Lemma 2.2.4}
If a spectrum $X$ is $H\F_p$-locally dualizable then
$L_{K(n)}X$ is $K(n)$-locally dualizable for each $0 < n < \infty$.
\endproclaim

\demo{Proof}
The Atiyah--Hirzebruch spectral sequence
$$
E^2_{s,t} = H_s(X; \pi_t K(n)) \Longrightarrow K(n)_{s+t}(X)
$$
shows that if $H_*(X; \F_p)$ is a (totally) finite $\F_p$-module,
then $K(n)_*(X)$ is a finitely generated $K(n)_*$-module for each $0 <
n < \infty$.
\qed
\enddemo

\remark{Remark 2.2.5}
More generally, for natural numbers $n < m$ it is not true that a
$K(m)$-locally dualizable spectrum $X$ must be $K(n)$-locally dualizable.
The lemma above would correspond to the case $m = \infty$.  The case
$X = K(n)$ provides an easy counterexample, with $K(m)_*(K(n)) = 0$
and $K(n)_*(K(n))$ infinitely generated over $K(n)_*$.  However, we are
principally interested in unreduced suspension spectra $X = \Sigma^\infty
T_+$, in which case the issue is: Does $K(m)_*(T)$ being finite in each
degree imply that $K(n)_*(T)$ is finite in each degree, for a topological
space $T$?  Replacing ``finite in each degree'' by ``trivial'' in this
statement, it becomes a theorem of Bousfield \cite{Bo99, 1.1}, with
a different proof under a finite type hypothesis by Wilson \cite{Wi99,
1.1}.  It is not clear to the author whether either of these proofs can
be adapted to resolve the stronger question.
\endremark

\subhead 2.3. Stably dualizable groups \endsubhead

Let $G$ be a topological group.  We write $S[G] = S \wedge G_+ = L_E
\Sigma^\infty G_+$ for the $E$-localization of the unreduced suspension
spectrum on $G$, and $DG_+ = F(S[G], S) = F(G_+, L_E \Sigma^\infty S^0)$
for its functional dual.  We may always suppose that $G$ is cofibrantly
based and of the homotopy type of a based CW-complex.

\definition{Definition 2.3.1}
A topological group $G$ is {\it $E$-locally stably dualizable\/}
if
$$
S[G] = L_E \Sigma^\infty G_+
$$
is dualizable in $\M_{S,E}$.
\enddefinition

\proclaim{Lemma 2.3.2}
The product $G = G_1 \times G_2$ of two $E$-locally stably dualizable
groups is again $E$-locally stably dualizable.
\endproclaim

\demo{Proof}
If $S[G_1]$ and $S[G_2]$ are dualizable, then so is $S[G] \cong S[G_1]
\wedge S[G_2]$, as remarked after Lemma~2.2.1.
\qed
\enddemo

The examples of Section~2.2 carry over as follows.  When $E$ is clear
from the context, we omit to say ``$E$-locally''.

\example{Example 2.3.3}
If $E = S$, then $G$ is a stably dualizable group if and only if $G_+$
is stably equivalent to a finite CW complex, up to an integer suspension.
So each compact Lie group $G$ is stably dualizable, since $G$ itself
is then a finite CW complex.
\endexample

\example{Example 2.3.4}
For $E = H\F_p$, a topological group $G$ is stably dualizable
if and only if $H_*(G; \F_p)$ is a (totally) finite $\F_p$-module.
\endexample

\example{Example 2.3.5}
For $E = K(n)$, a topological group $G$ is stably dualizable if
and only if $K(n)_*(G)$ is a finitely generated $K(n)_*$-module.
By the calculations of Ravenel and Wilson \cite{RaW80, 11.1} for $p$
odd, and \cite{JW85, Appendix} for $p=2$, each Eilenberg--Mac\,Lane
space $G = K(\pi, q) = B^q \pi$ for $\pi$ a finite abelian group is
a stably dualizable group.  More generally, by \cite{HRW98, 1.1} any
topological group $G$ with only finitely many nonzero homotopy groups,
each of which is a finite abelian $p$-group, has finite $K(n)$-homology,
hence is stably dualizable.
\endexample

\remark{Remark 2.3.6}
By Lemma~2.2.4, compact Lie groups or $p$-compact groups provide
examples of $K(n)$-locally stably dualizable groups, since they have
finite mod~$p$ homology, and therefore have finite Morava $K$-theory.
The Eilenberg--Mac\,Lane space examples above, for $q\ge1$, do not
arise in this fashion, since they have infinite mod~$p$ homology.
By Example~2.3.4 they do not extend to stably dualizable groups in the
$p$-complete or integral category.

Building on Remark~2.2.5, if it turns out that $K(m)_*(G)$ being
finite over $K(m)_*$ implies that $K(n)_*(G)$ is finite over $K(n)_*$,
for topological groups $G$ and natural numbers $n < m$, then each
$K(m)$-locally stably dualizable group will also be a $K(n)$-locally
stably dualizable group.  By \cite{Bo99, 1.1}, each $K(m)$-equivalence
$G_1 \to G_2$ is also a $K(n)$-equivalence, so there will then be a
``chromatic'' tower of $K(n)$-equivalence classes of $K(n)$-locally
stably dualizable groups, for $1 \le n \le \infty$, with maps
from the set at height $m$ to the set at height $n$, for all $n < m$.
\endremark

\subhead 2.4. $E$-compact groups \endsubhead

The material in this section is included to enable a more precise
comparison with the Dwyer--Wilkerson theory of $p$-compact groups,
but is not needed elsewhere in the paper.

Suppose that the spectrum $E$ is in fact an $S$-algebra \cite{EKMM97}.
This includes all the examples $E = HR$ for rings $R$, $E = S$, $E =
E(n)$ and $E = K(n)$ considered above, although the $S$-algebra structure
on e.g.~$K(n)$ is not at all unique.  We now consider a version of
Bousfield--Kan $p$-completion for the homology theory represented by
$E$, following \cite{BCM78} and \cite{BT00, \S2}.  Let $\Omega^\infty E$
denote the underlying infinite loop space of $E$, so that $\Sigma^\infty$
is left adjoint to $\Omega^\infty$.  The following terminology extends
that of \cite{BK72, \S I.5}.

\definition{Definition 2.4.1}
(a)
Let $E(X) = \Omega^\infty(E \wedge \Sigma^\infty X)$ define an endofunctor
of based topological spaces, with $\pi_* E(X) = \widetilde E_*(X)$ for
all $*\ge0$.  The $S$-algebra multiplication $\mu \: E \wedge E \to E$
and the adjunction counit $\Sigma^\infty \Omega^\infty E \to E$ induce
a natural transformation $\mu_X \: (E \circ E)(X) = E(E(X)) \to E(X)$.
The $S$-algebra unit $\eta \: S \to E$ and the adjunction unit $X \to
Q(X) = \Omega^\infty \Sigma^\infty X$ induce a natural transformation
$\eta_X \: X \to E(X)$.  These make $(E(-), \mu, \eta)$ a monad.

(b)
For each space $X$, the {\it $E$-completion} $X^\wedge_E = \Tot
E(X)^\bullet$ is defined as the totalization of the cosimplicial space
$$
[q] \mapsto E(X)^q = (E \circ \dots \circ E)(X)
$$
(the $(q+1)$-fold iterate of the endofunctor $E(-)$), with coface and
codegeneracy maps induced by $\mu$ and $\eta$, respectively.  There is a
natural map $X \to X^\wedge_E$.  We say that $X$ is {\it $E$-complete\/}
if this map is a weak equivalence.

(c)
Each space $E(X)^q$ is $E$-local, so $X^\wedge_E$ is $E$-local and
there is a canonical factorization $X \to L_E X \to X^\wedge_E$ through
the Bousfield localization to the completion.  We say that $X$ is {\it
$E$-good\/} if $X \to X^\wedge_E$ is an $E$-equivalence, or equivalently,
if $L_E X \simeq X^\wedge_E$.
\enddefinition

\definition{Definition 2.4.2}
An {\it $E$-compact group\/} is an $E$-locally stably dualizable group
$G$ whose classifying space $BG \simeq (BG)^\wedge_E$ is $E$-complete.
\enddefinition

\example{Example 2.4.3}
If $E = S$, then $S$-completion equals $H\Z$-completion \cite{Ca91,
II.3}, so nilpotent spaces are $S$-complete.  A connected compact Lie
group $G$ has simply-connected, thus nilpotent, classifying space $BG$,
so such a group $G$ is also an $S$-compact group.
\endexample

\example{Example 2.4.4}
If $E = H\F_p$, $H\F_p$-completion equals $p$-completion, so a topological
group $G$ is $H\F_p$-compact if and only if $G \simeq \Omega BG$ is a
$p$-compact group in the sense of Dwyer and Wilkerson \cite{DW94}.

The loop space $\Omega(BG)^\wedge_p$ of the $p$-completed classifying
space of a connected compact Lie group provides some examples of a
$p$-compact group, but there are also exotic examples, such as (i)
the $p$-compact Sullivan sphere
$$
(S^{2p-3})^\wedge_p = \Omega(EW \times_W BA)^\wedge_p
$$
for $p$ odd, with $A = B\Z_p$, $BA = K(\Z_p, 2)$ and $W = (\Z/p)^*$ acting
on $A$ through multiplication by the $p$-adic roots of unity, and (ii)
the $2$-compact Dwyer--Wilkerson group $DI(4)$ \cite{DW93}.  With the
exception of $(S^3)^\wedge_3$ these only exist locally, in the sense
that they do not extend to integrally defined stably dualizable groups.
\endexample

\example{Example 2.4.5}
For $q \ge n$ and $\pi$ a finite abelian group the Eilenberg--Mac\,Lane
spaces $G = K(\pi, q)$ have $K(n)_*(BG) = K(n)_*$ by \cite{RaW80, 11.1},
hence these can never be nontrivial $K(n)$-compact groups.  When $0 \le
q < n$ and $\pi$ is a finite abelian $p$-group it is known that $BG =
K(\pi, q+1)$ is $K(n)$-local by \cite{Bo82, 7.4}, so if these classifying
spaces are also $K(n)$-good, then the Eilenberg--Mac\,Lane spaces $K(\pi,
q)$ are $K(n)$-compact groups.
\endexample

\remark{Remark 2.4.6}
It is to be expected that connected compact Lie groups, $p$-compact
groups or the more exotic $K(n)$-locally stably dualizable groups
of Example~2.3.5 provide examples of $K(n)$-compact groups through
$K(n)$-completion at the level of classifying spaces.  However, it is
not clear in what generality the natural map
$$
G \to \Omega (BG)^\wedge_{K(n)}
$$
is a $K(n)$-equivalence.  There is a $K(n)$-based
Eilenberg--Moore spectral sequence \cite{JO99}
$$
E_2^{s,t} = \widehat{\Tor} {}^{s,t}_{K(n)^*(Y)}(K(n)^*, K(n)^*)
\Longrightarrow K(n)^{s+t}(\Omega Y) \,,
\tag 2.4.7
$$
where $\widehat{\Tor} {}^s = \widehat{\Tor}_{-s}$ are the left derived
functors of the completed tensor product, but too little is known
about its convergence.  Certainly the space $Y$ should be $K(n)$-local
for this to have a chance, but by analogy with the mod~$p$ case,
it is more plausible that the correct condition is that $Y$ should
be $K(n)$-complete.  Since a $K(n)$-complete space is a limit of
spaces of the form $\Omega^\infty(K(n) \wedge T)$, it may suffice to
verify convergence for such spaces, or for the individual spaces $Y =
\underline{K(n)}_q$ of the $\Omega$-spectrum of $K(n)$, for $q\ge0$.

Bauer has made some progress in this direction.  So if (i) $G$ is a
$K(n)$-locally stably dualizable group, (ii) $BG$ is $K(n)$-good, (iii)
$K(n)^*(BG)$ is a finitely generated power series ring over $K(n)^*$,
or more generally, $\widehat{\Tor} {}^{**}_{K(n)^*(BG)}(K(n)^*, K(n)^*)$
is finite over $K(n)^*$, and (iv) the $K(n)$-based Eilenberg--Moore
spectral sequence for $Y = (BG)^\wedge_{K(n)}$ converges, then it follows
that $\Omega Y = \Omega(BG)^\wedge_{K(n)}$ is a $K(n)$-compact group.
\endremark

\remark{Remark 2.4.8}
The examples that are abelian topological groups can be expected to play
a similar role to that of tori in the theory of compact Lie groups.
For non-abelian examples it is natural to look to finite Postnikov
systems, as in \cite{HRW98}, or to looped completed Borel constructions
of the form
$$
G = \Omega (EW \times_W BA)^\wedge_{K(n)}
$$
where $A$ is an abelian topological group, such as $A = K(\pi, q)$,
the Weyl group $W$ is a finite group acting on $A$ from the left,
$EW \times_W BA = B(W \ltimes A)$ is the classifying space of the
semi-direct product $W \ltimes A$ and $(-)^\wedge_{K(n)}$ denotes the
$K(n)$-completion of spaces.  To analyze the $K(n)$-homology of $G$ it is
again necessary to study the convergence properties of the $K(n)$-based
Eilenberg--Moore spectral sequence in the path--loop fibration of $Y =
B(W \ltimes A)^\wedge_{K(n)}$.
\endremark

\remark{Remark 2.4.9}
In his Master's thesis, H{\aa}kon Schad Bergsaker \cite{Be05, 6.6}
has shown that for $A = K(\Z_p, n)$ and $W = (\Z/p)^*$ acting through
multiplication by the $p$-adic roots of unity, the construction
above produces a loop space (or topological group) $G$ with the Morava
$K$-theory of a $(2m-1)$-sphere, for $m = (p^n-1)/\gcd(p-1,n)$, under the
assumption that the appropriate $K(n)$-based Eilenberg--Moore spectral
sequence~(2.4.7) converges.  Conversely, he shows \cite{Be05, 4.15}
that for each prime~$p$ and height $n$ there are only finitely many
$m$ for which $G = L_{K(n)} S^{2m-1}$ can be a $K(n)$-locally stably
dualizable group, assuming the existence of a map $t \: B\Z/p \to BG$
with nontrivial $K(n)^*(t)$.
\endremark

\subhead 2.5. The dualizing and inverse dualizing spectra \endsubhead

Let $EG = B(*, G, G)$ be the usual free, contractible right $G$-space.
Let $X$ be a spectrum with right $G$-action, and let $Y$ be a spectrum
with left $G$-action.  We define the {\it $G$-homotopy fixed points\/}
of $X$ to be
$$
X^{hG} = F(EG_+, X)^G
$$
and the {\it $G$-homotopy orbits\/} of $Y$ to be
$$
Y_{hG} = EG_+ \wedge_G Y \,.
$$
In all cases $G$ acts on $EG$ from the right.  These constructions only
involve naive $G$-equivariant spectra, or spectra with $G$-action, in the
sense that no deloopings with respect to non-trivial $G$-representations
are involved.  Each $G$-equivariant map $X_1 \to X_2$ that is an
equivalence induces an equivalence $X_1^{hG} \to X_2^{hG}$ of homotopy
fixed points, and similarly for homotopy orbits.

\definition{Definition 2.5.1}
Let $G$ be an $E$-locally stably dualizable group.  The group
multiplication provides the suspension spectrum $S[G] = L_E \Sigma^\infty
G_+$ with mutually commuting {\it standard\/} left and right $G$-actions.
We define the {\it dualizing spectrum} $S^{adG}$ of $G$ to be the
$G$-homotopy fixed point spectrum
$$
S^{adG} = S[G]^{hG} = F(EG_+, S[G])^G
$$
of $S[G]$, formed with respect to the standard right $G$-action.
The standard left action on $S[G]$ induces a left $G$-action on $S^{adG}$.
\enddefinition

\remark{Remark 2.5.2}
A discrete group $G$ of type $FP$ (e.g.~the classifying space $BG$ is
finitely dominated) is called a {\it duality group\/} if $H^*(G; \Z[G])$
is concentrated in a single degree~$n$ and torsion free.  The $G$-module
$D = H^n(G; \Z[G])$ is then called the {\it dualizing module\/} of $G$,
cf.~\cite{Brn82, VIII.10}.   The spectrum level construction above is
analogous to this algebraic notion, and was previously considered for
topological groups by Dwyer and by Klein \cite{Kl01, \S1}, and for
$p$-compact groups by Bauer \cite{Ba04, 4.1}.  However, in algebra the
focus is on $G$ with a finiteness condition on $BG$, whereas in the
topological cases $G$ itself satisfies a finiteness condition.  Klein
writes $D_G$ and Bauer writes $S_G$ for the dualizing spectrum of $G$.
We use $D$ for the functional dual and $S$ for the sphere spectrum, so
we prefer to write $S^{adG}$ instead, in view of the compact Lie group
example recalled immediately below.  Our construction differs a tiny
bit from that of Dwyer and Klein, due to our implicit $E$-localization.
\endremark

\example{Examples 2.5.3}
(a)
When $G$ is a finite group, there is a canonical equivalence $S[G] =
S \wedge G_+ \simeq F(G_+, S)$, so $S[G]^{hG} \simeq F(G_+, S)^{hG}
\cong F(EG_+, S) \simeq S$ is naturally equivalent to the sphere spectrum.

(b)
More generally, when $G$ is a compact Lie group Klein \cite{Kl01, 10.1}
shows that the dualizing spectrum $S^{adG}$ is equivalent as a spectrum
with left $G$-action to the suspension spectrum of the representation
sphere associated to the adjoint representation $adG$ of $G$, i.e.,
the left conjugation action of $G$ on its tangent space $T_eG$ at the
identity.

(c)
In the case of a $p$-compact group $G$, Bauer \cite{Ba04} shows that
$S^{adG} \simeq (S^d)^\wedge_p$ for some integer $d = \dim_p G$ called
the {\it $p$-dimension\/} of $G$, and that $S^{adG}$ takes over the
role of the representation sphere in the duality theory in that context.
The present paper extends some of Bauer's work to the $E$-local stable
category.
\endexample

\proclaim{Lemma 2.5.4}
When $G$ is abelian, the left $G$-action on $S^{adG}$ is homotopically
trivial, in the sense that it extends over the inclusion $G \subset EG$
to an action by the contractible topological group $EG$.
\endproclaim

\demo{Proof}
When $G$ is abelian, the left and right $G$-actions on $S[G]$ agree.
In $S^{adG} = F(EG_+, S[G])^G$ the right action on $S[G]$ is
equal to the right action on $EG_+$, which in the commutative case
factors as
$$
EG_+ \wedge G_+ \subset EG_+ \wedge EG_+ @>>> EG_+ \,.
\qed
$$
\enddemo

\remark{Remark 2.5.5}
It can be inconvenient to study the $E$-homology of $S^{adG}$ directly
from its definition as a homotopy fixed point spectrum.  We shall soon
see that this dualizing spectrum is the functional dual of another
spectrum $S^{-adG}$, which we call the inverse dualizing spectrum,
and which admits a computationally more convenient construction as
a homotopy orbit spectrum.  Once we know that these two spectra are
indeed dualizable, and mutually dual, this provides a convenient route
to $E$-homological calculations.
\endremark

\definition{Definition 2.5.6}
Let $G$ be a stably dualizable group.  The left and right $G$-actions
on $S[G]$ induce {\it standard\/} right and left $G$-actions on its
functional dual $DG_+ = F(S[G], S)$, respectively, by acting in the source
of the mapping spectrum.  We define the {\it inverse dualizing spectrum}
$S^{-adG}$ of $G$ to be the $G$-homotopy orbit spectrum
$$
S^{-adG} = (DG_+)_{hG} = EG_+ \wedge_G DG_+
$$
of $DG_+$, formed with respect to the standard left $G$-action.  These
left and right actions commute, so the standard right action on $DG_+$
induces a right $G$-action on $S^{-adG}$.
\enddefinition

\proclaim{Proposition 2.5.7}
There is a natural equivalence
$$
S^{adG} \simeq DS^{-adG}
$$
between the dualizing spectrum and the functional dual of the inverse
dualizing spectrum, as spectra with left $G$-action.
\endproclaim

\demo{Proof}
The canonical equivalence $\rho \: S[G] \to DDG_+ = F(DG_+, S)$
induces an equivalence $\rho^{hG}$ of $G$-homotopy fixed points,
from $S^{adG}$ to
$$
F(DG_+, S)^{hG} = F(EG_+, F(DG_+, S))^G
\cong F(EG_+ \wedge_G DG_+, S) = DS^{-adG} \,.
\qed
$$
\enddemo

\head 3. Duality theory \endhead

\subhead 3.1. Poincar{\'e} duality \endsubhead

Let $G$ be a stably dualizable group.  The topological group structure on
$G$ makes $S[G]$ a cocommutative Hopf $S$--algebra, with product $\phi
\: S[G] \wedge S[G] \to S[G]$, unit $\eta \: S \to S[G]$, coproduct
$\psi \: S[G] \to S[G] \wedge S[G]$, counit $\epsilon \: S[G] \to S$
and conjugation (antipode) $\chi \: S[G] \to S[G]$, induced by the group
multiplication $m \: G \times G \to G$, unit inclusion $\{e\} \to G$,
diagonal map $\Delta \: G \to G \times G$, collapse map $G \to \{e\}$
and group inverse $i \: G \to G$, respectively.

The product $\phi$ and unit $\eta$ makes $S[G]$ an $E$-local $S$-algebra
in $\M_{S,E}$, while the coproduct, counit and conjugation need only be
defined in the $E$-local stable category $\D_{S,E}$.

The standard right $G$-action on $DG_+$ makes $DG_+$ a right
$S[G]$-module.  The module action is given by the map
$$
\alpha \: DG_+ \wedge S[G] @>>> DG_+
$$
that in symbols takes $\xi \wedge x$ to $\xi * x \: y \mapsto \xi(xy)$.
Inspired by \cite{Ba04, \S4.3}, we consider the following shearing
equivalence.  Its definition is simpler than that considered by Bauer,
but the key idea is the same.

\definition{Definition 3.1.1}
Let the {\it shear map} $sh \: DG_+ \wedge S[G] \to DG_+ \wedge
S[G]$ be the composite map
$$
sh \: DG_+ \wedge S[G] @>1\wedge\psi>> DG_+ \wedge S[G] \wedge
S[G] @>\alpha\wedge1>> DG_+ \wedge S[G] \,.
$$
Algebraically, $sh \: \xi \wedge x \mapsto \sum (\xi * x') \wedge x''$
where $\psi(x) = \sum x' \wedge x''$.
\enddefinition

The standard left and right $G$-actions on $S[G]$ (and $DG_+$) can
be converted into right and left $G$-actions on $S[G]$ (and $DG_+$),
respectively, by way of the group inverse $i \: G \to G$.  We refer
to these non-standard actions as {\it actions through inverses}.
For example, the left $G$-action through inverses on $DG_+$ is given by
the composite map
$$
S[G] \wedge DG_+ @>\gamma>\cong> DG_+ \wedge S[G] @>1\wedge\chi>> DG_+
\wedge S[G] @>\alpha>> DG_+ \,,
$$
where $\gamma \: X \wedge Y \to Y \wedge X$ denotes the canonical twist
map.  Algebraically, this action takes $(x, \xi)$ to $\xi * \chi(x) \:
y \mapsto \xi(\chi(x) y)$.

\proclaim{Lemma 3.1.2}
The shear map $sh$ is equivariant with respect to each of
the following three mutually commuting $G$-actions:
\roster
\item
The first, left $G$-action given by the action through inverses on
$DG_+$ and the standard action on $S[G]$ in the source, and the
standard action on $S[G]$ in the target;
\item
The second, right $G$-action given by the action through inverses
on $DG_+$ in the source, and the action through inverses on
$DG_+$ in the target;
\item
The third, right $G$-action given by the standard action on $S[G]$
in the source and by the standard actions on $DG_+$ and $S[G]$
in the target.
\endroster
Each action is trivial on the remaining smash factors.
\endproclaim

\demo{Proof}
In each case this is clear by inspection.
\qed
\enddemo

\proclaim{Lemma 3.1.3}
The shear map $sh$ is an equivalence, with homotopy inverse
given by the composite map
$$
DG_+ \wedge S[G] @>1\wedge\psi>> DG_+ \wedge S[G] \wedge S[G]
@>1\wedge\chi\wedge1>> DG_+ \wedge S[G] \wedge S[G] @>\alpha\wedge1>>
DG_+ \wedge S[G] \,.
$$
\endproclaim

\demo{Proof}
This is an easy diagram chase, using coassociativity of $\psi$, the fact
that $\alpha$ is a right $S[G]$-module action with respect to the product
$\phi$ on $S[G]$, the Hopf conjugation identities $\phi(\chi\wedge1)\psi
\simeq \eta\epsilon \simeq \phi(1\wedge\chi)\psi$, counitality for $\psi$
and unitality for $\alpha$.
\qed
\enddemo

\proclaim{Theorem 3.1.4}
Let $G$ be a stably dualizable group.
There is a natural equivalence
$$
DG_+ \wedge S^{adG} @>\simeq>> S[G] \,.
$$
It is equivariant with respect to the first, left $G$-action through
inverses on $DG_+$, the standard left action on $S^{adG}$ and the
standard left action on $S[G]$.  It is also equivariant with respect
to the second, right $G$-action through inverses on $DG_+$, the trivial
action on $S^{adG}$ and the standard right action on $S[G]$.
\endproclaim

\demo{Proof}
The shear equivalence $sh \: DG_+ \wedge S[G] \to DG_+ \wedge S[G]$
induces a natural equivalence
$$
(sh)^{hG} \: (DG_+ \wedge S[G])^{hG} @>\simeq>> (DG_+ \wedge S[G])^{hG}
$$
of $G$-homotopy fixed points with respect to the third, right $G$-action.
Note that this action is different in the source and in the target
of $sh$.

There is a natural equivalence to the source of $(sh)^{hG}$:
$$
DG_+ \wedge S^{adG} = DG_+ \wedge S[G]^{hG} @>\simeq>> (DG_+ \wedge
S[G])^{hG} \,.
$$
To see that this map is an equivalence, consider the commutative square
$$
\xymatrix{
DG_+ \wedge S[G]^{hG} \ar[r] \ar[d]^{\simeq} &
(DG_+ \wedge S[G])^{hG} \ar[d]^{\simeq} \\
F(G_+, S[G]^{hG}) \ar[r]^{\cong} &
F(G_+, S[G])^{hG} \,.
}
$$
The vertical maps are equivalences, because $S[G]$ is dualizable and
passage to homotopy fixed points respects equivalences.  Hence the upper
horizontal map is also an equivalence.

There is also a (composite) natural equivalence from the target of
$(sh)^{hG}$:
$$
(DG_+ \wedge S[G])^{hG} @>\simeq>> F(G_+, S[G])^{hG}
@>\simeq>> S[G] \,.
$$
The left hand map is an equivalence because $S[G]$ is dualizable, by the
same argument as above.  The right hand map is the composite equivalence
$$
F(G_+, S[G])^{hG} \cong F(EG_+ \wedge G_+, S[G])^G
\cong F(EG_+, S[G]) @>\simeq>> S[G] \,.
$$
Here the middle isomorphism uses that $G$ acts freely on $G_+$ in
the source.

The composite of these three natural equivalences is the desired natural
equivalence $DG_+ \wedge S^{adG} \to S[G]$.  The equivariance statements
follow by inspection.
\qed
\enddemo

\remark{Remark 3.1.5}
We call $DG_+ \wedge S^{adG} \simeq S[G]$ the {\it Poincar{\'e} duality
equivalence}.  It shows how $S[G]$ is functionally self-dual, up to a
shift by the dualizing spectrum.  See also Remark~3.3.5.  The equivariance
statements in the theorem express the standard left and trivial right
$G$-actions on $S^{adG}$ in terms of the more familiar $G$-actions on
$DG_+$ and $S[G]$.
\endremark

\proclaim{Lemma 3.1.6}
Let $G_1$ and $G_2$ be stably dualizable groups.  There is a natural
equivalence
$$
S^{ad G_1} \wedge S^{ad G_2} \simeq S^{ad(G_1 \times G_2)}
$$
of spectra with standard left (and trivial right) $(G_1 \times
G_2)$-actions.
\endproclaim

\demo{Proof}
The Poincar{\'e} duality equivalences for $G_1$, $G_2$ and
$(G_1 \times G_2)$ compose to an equivalence
$$
\multline
DG_{1+} \wedge S^{ad G_1} \wedge DG_{2+} \wedge S^{ad G_2}
\simeq S[G_1] \wedge S[G_2] \\ \simeq S[G_1 \times G_2]
\simeq D(G_1 \times G_2)_+ \wedge S^{ad (G_1 \times G_2)}
\,.
\endmultline
$$
It is equivariant with respect to the first, left $(G_1 \times
G_2)$-action that involves the standard left action on $S^{ad G_1}$,
$S^{ad G_2}$ and $S^{ad (G_1 \times G_2)}$, as well as with respect to
the second, right $(G_1 \times G_2)$-action through inverses on $DG_{1+}
\wedge DG_{2+}$ and $D(G_1 \times G_2)_+$.  Taking homotopy fixed
points with respect to the second, right action we obtain the desired
equivalence, which is equivariant with respect to the first, left action.
\qed
\enddemo

\remark{Remark 3.1.7}
A similar relation $S^{adG} \simeq S^{adH} \wedge S^{adQ}$ is likely to
hold for an extension $1 \to H \to G \to Q \to 1$ of stably dualizable
groups, cf.~\cite{Kl01, Thm.~C}, but for simplicity we omit the then
necessary discussion of how to promote $S^{ad H}$ to a spectrum with
$G$-action, etc.
\endremark

\subhead 3.2. Inverse Poincar{\'e} duality \endsubhead

The aim of this section is to establish an inverse Poincar{\'e}
equivalence
$$
S[G] \wedge S^{-adG} \simeq DG_+ \,.
$$
The initial idea is to functionally dualize the construction of the shear
map in Section~3.1, and to apply homotopy orbits in place of homotopy
fixed points.  Following Milnor and Moore \cite{MM65, \S3}, we identify
the functional dual of a smash product $X \wedge Y$ of dualizable spectra
with the smash product $DX \wedge DY$, in that order, via the canonical
equivalence
$$
DX \wedge DY = F(X, S) \wedge F(Y, S) @>\wedge>\simeq>
F(X \wedge Y, S \wedge S) = D(X \wedge Y) \,.
$$
However, to form homotopy orbits we need genuine $G$-equivariant maps,
and it is generally not the case that a $G$-equivariant inverse can be
found for the (weak) equivalence displayed above.  Thus some care will
be in order.

Working for a moment in the $E$-local stable category $\D_{S,E}\ =
\bar h \M_{S,E}$, let
$$
\beta \: S[G] @>>> S[G] \wedge DG_+
$$
be dual to the module action map $\alpha \: DG_+ \wedge S[G] \to DG_+$.
It makes $S[G]$ a right $DG_+$-comodule spectrum, up to homotopy, where
$DG_+$ has the weakly defined coproduct $\psi' \: DG_+ \to DG_+ \wedge
DG_+$ that is dual to $\phi$.  Furthermore, let
$$
\phi' \: DG_+ \wedge DG_+ @>>> DG_+
$$
be the (strictly defined) product on $DG_+$ that is dual to $\psi$.
The functional dual $sh^\#$ of the shear map is then the composite
$$
sh^\# \: S[G] \wedge DG_+ @>\beta\wedge1>> S[G] \wedge DG_+ \wedge DG_+
@>1\wedge\phi'>> S[G] \wedge DG_+ \,,
$$
which is an equivalence by Lemma~3.1.3 and duality.

Returning to the category $\M_{S,E}$, we shall now obtain $G$-equivariant
representatives for these maps.

\definition{Definition 3.2.1}
Let $\tilde\phi \: S[G] \to F(S[G], S[G])$ be right adjoint to the
opposite product map $\phi \gamma \: S[G] \wedge S[G] \to S[G]$.
Algebraically, $\tilde\phi \: x \mapsto (y \mapsto yx)$.  Let $\psi^\#
\: F(S[G] \wedge S[G], S[G] \wedge S) \to F(S[G], S[G])$ be given by
precomposition by $\psi \: S[G] \to S[G] \wedge S[G]$ and postcomposition
with $S[G] \wedge S \cong S[G]$.

The {\it dual shear map} $sh' \: S[G] \wedge DG_+ \to F(S[G],
S[G])$ is defined to be the composite map:
$$
\multline
sh' \: S[G] \wedge DG_+ @>\tilde\phi\wedge1>> F(S[G], S[G]) \wedge DG_+ \\
@>\wedge>\simeq> F(S[G] \wedge S[G], S[G] \wedge S) @>\psi^\#>> F(S[G], S[G])
\,.
\endmultline
$$
It is equivariant with respect to the left $G$-action given by the
standard left actions on $S[G]$ and $DG_+$ on the left hand side, and
the left action through the standard right action on the $S[G]$ in the
source of the mapping spectrum.
\enddefinition

\proclaim{Theorem 3.2.2}
The dual shear map $sh'$ is homotopic to the composite map
$$
S[G] \wedge DG_+ @>sh^\#>\simeq> S[G] \wedge DG_+ @>\nu\gamma>\simeq>
F(S[G], S[G]) \,.
$$
In particular, $sh'$ is an equivalence.  On $G$-homotopy orbit spectra
it induces an equivalence
$$
DG_+ \simeq S[G] \wedge S^{-adG} \,.
$$
\endproclaim

\demo{Proof}
The right action map $\alpha$ factors up to homotopy as
the composite
$$
\multline
DG_+ \wedge S[G]
@>\psi'\wedge1>> DG_+ \wedge DG_+ \wedge S[G] \\
@>1\wedge\gamma>> DG_+ \wedge S[G] \wedge DG_+
@>\epsilon\wedge1>> S \wedge DG_+ = DG_+ \,.
\endmultline
$$
Here $\epsilon \: DG_+ \wedge S[G] \to S$ is the pairing that evaluates
a function on an element in its source.  Let $\eta \: S[G] \wedge DG_+
\to S$ be its functional dual, in the homotopy category.  Then the dual
map $\beta$ factors up to homotopy as
$$
\multline
S[G] \cong S \wedge S[G] @>\eta\wedge1>> S[G] \wedge DG_+ \wedge S[G] \\
@>1\wedge\gamma>> S[G] \wedge S[G] \wedge DG_+
@>\phi\wedge1>> S[G] \wedge DG_+ \,.
\endmultline
$$
A diagram chase then verifies that $\tilde\phi$ is homotopic
to the composite
$$
S[G] @>\beta>> S[G] \wedge DG_+ @>\gamma>\cong>
DG_+ \wedge S[G] @>\nu>\simeq> F(S[G], S[G]) \,.
$$
A similar chase shows that the diagram
$$
\xymatrix{
S[G] \wedge DG_+ \wedge DG_+ \ar[rr]^-{1\wedge\phi'}
\ar[d]^-{\nu\gamma\wedge1}_-{\simeq} &&
S[G] \wedge DG_+ \ar[d]^-{\nu\gamma}_-{\simeq} \\
F(S[G], S[G]) \wedge DG_+ \ar[r]^-{\wedge}_-{\simeq} &
F(S[G] \wedge S[G], S[G] \wedge S) \ar[r]^-{\psi^\#} &
F(S[G], S[G])
}
$$
homotopy commutes.

Taken together, these diagrams show that $\nu\gamma \circ sh^\# \simeq
sh'$.  Applying $G$-homotopy orbits to the chain of equivalences
$$
S[G] \wedge DG_+ @>sh'>\simeq> F(S[G], S[G]) @<\nu\gamma<\simeq<
S[G] \wedge DG_+
$$
we obtain the desired chain of equivalences
$$
\multline
DG_+ \simeq (S[G] \wedge DG_+)_{hG}
@>(sh')_{hG}>\simeq> F(S[G], S[G])_{hG} \\
@<(\nu\gamma)_{hG}<\simeq< (S[G] \wedge DG_+)_{hG}
\simeq S[G] \wedge S^{-adG} \,.
\qed
\endmultline
$$
\enddemo

\proclaim{Proposition 3.2.3}
Let $G$ be a stably dualizable group.  The dualizing spectrum
$S^{adG}$ and the inverse dualizing spectrum $S^{-adG}$ are both
dualizable spectra.  Hence
$$
S^{-adG} \simeq DS^{adG}
$$
as spectra with right $G$-action.  The inverse Poincar{\'e} equivalence
$$
S[G] \wedge S^{-adG} \simeq DG_+
$$
is equivariant with respect to the dual $G$-actions to those of
Theorem~3.1.4:  The first of these is the right $G$-action through
inverses on $S[G]$, the standard right action on $S^{-adG}$ and the
standard right action on $DG_+$.  The second is the left $G$-action
through inverses on $S[G]$, the trivial action on $S^{-adG}$ and the
standard left action on $DG_+$.
\endproclaim

\demo{Proof}
It suffices to prove that $S^{-adG}$ is dualizable, in view of
Proposition~2.5.7 and Theorem~3.1.4.  We must show that the canonical
map
$$
\nu \: DS^{-adG} \wedge S^{-adG} @>>> F(S^{-adG}, S^{-adG})
$$
is an equivalence.  We first check that $\nu$ smashed with the
identity map of $S[G]$ is an equivalence.  This map factors as
the composite
$$
\multline
DS^{-adG} \wedge S^{-adG} \wedge S[G]
\simeq DS^{-adG} \wedge DG_+
@>\nu>\simeq> F(S^{-adG}, DG_+) \\
\simeq F(S^{-adG}, S^{-adG} \wedge S[G])
@<\nu<\simeq< F(S^{-adG}, S^{-adG}) \wedge S[G] \,.
\endmultline
$$
Here the first and third equivalences follow from the inverse Poincar{\'e}
equivalence, while the second and fourth equivalences are consequences of
the dualizability of $DG_+$ and $S[G]$, respectively.  Thus $\nu \wedge
1_{S[G]}$ is an equivalence.  Since $S$ is a retract of $S[G]$, it
follows that also $\nu$ itself is an equivalence.  Hence $S^{-adG}$
is dualizable.
\qed
\enddemo

\subhead 3.3. The Picard group \endsubhead

The Picard group of the category of $E$-local $S$-modules was introduced
by M.~Hopkins; see \cite{HMS94}.

\definition{Definition 3.3.1}
An $E$-local $S$-module $X$ is {\it invertible\/} if there exists a
spectrum $Y$ with $X \wedge Y \simeq S$ in $\M_{S,E}$.  Then $Y$ is
also invertible.  The smash product $X \wedge X'$ of two invertible
spectra $X$ and $X'$ is again invertible.

The {\it $E$-local Picard group} $\Pic_E = \Pic(\M_{S,E})$ is the set
of equivalence classes of invertible $E$-local $S$-modules.  We write
$[X] \in \Pic_E$ for the equivalence class of $X$.  The abelian group
structure on $\Pic_E$ is defined by $[X] + [X'] = [X \wedge X']$
and $-[X] = [Y]$, with $X$, $Y$ and $X'$ as above.
\enddefinition

\example{Example 3.3.2}
The only invertible spectra in $\M_S$ are the sphere spectra $S^d =
\Sigma^d S$ for integers $d \in \Z$, so $\Pic_S \cong \Z$.  Similarly,
in the $p$-complete category $\M_{S, H\F_p}$ the invertible spectra are
precisely the $p$-completed sphere spectra $(S^d)^\wedge_p$ for
$d \in \Z$, so $\Pic_{H\F_p} \cong \Z$ too.
\endexample

\example{Example 3.3.3}
By Hopkins, Mahowald and Sadofsky \cite{HMS94, 1.3}, a $K(n)$-local
spectrum $X$ is invertible if and only if $K(n)_*(X)$ is free of rank
one over $K(n)_*$.  These authors show \cite{HMS94, 2.1, 2.7, 3.3}
that for $n=1$ and $p\ne2$ there is a non-split extension
$$
0 \to \Z_p^\times @>>> \Pic_{K(1)} @>>> \Z/2 \to 0
$$
while for $n=1$ and $p=2$ there is a non-split extension
$$
0 \to \Z_2^\times @>>> \Pic_{K(1)} @>>> \Z/8 \to 0 \,.
$$
Furthermore, they show \cite{HMS94, 7.5} that when $n^2 \le 2p-2$ and
$p > 2$ there is an injection $\alpha \: \Pic_{K(n)} \to H^1(\Bbb S_n;
\pi_0(E_n)^\times)$, where $E_n$ is the Hopkins--Miller commutative
$S$-algebra and $\Bbb S_n$ is (one of the variants of) the $n$-th
Morava stabilizer group.  This permits an algebraic identification of
$\Pic_{K(2)}$ for $p$ odd.  The homomorphism $\alpha$ seems to have a
non-trivial kernel for $n=2$ and $p=2$, cf.~\cite{HMS94, \S6}.
\endexample

\proclaim{Theorem 3.3.4}
Let $G$ be a stably dualizable group.  Then
$$
S^{adG} \wedge S^{-adG} \simeq S
$$
so $S^{adG}$ and $S^{-adG}$ are mutually inverse invertible spectra in
the $E$-local stable category.  Hence the equivalence classes $[S^{adG}]$
and $[S^{-adG}]$ represent inverse elements in the $E$-local Picard
group $\Pic_E$.
\endproclaim

\demo{Proof}
The Poincar{\'e} duality equivalence and the inverse Poincar{\'e}
equivalence provide a chain of equivalences
$$
S[G] \wedge S^{-adG} \wedge S^{adG} \simeq DG_+ \wedge S^{-adG}
\simeq S[G] \,,
$$
which is equivariant with respect to the standard left action on both
copies of $S[G]$, the trivial action on $S^{-adG}$ and the standard left
action on $S^{adG}$.  Taking $G$-homotopy orbits of both sides yields
the required equivalence
$$
S^{-adG} \wedge S^{adG} \simeq S[G]_{hG} \simeq S \,.
\qed
$$
\enddemo

\remark{Remark 3.3.5}
These results show that the shift given by smashing with $S^{adG}$,
as in the Poincar{\'e} duality equivalence, is really an invertible
self-equivalence of the stable homotopy category of spectra with
$G$-action, in that it can be undone by smashing with $S^{-adG} \simeq
DS^{adG}$.
\endremark

\definition{Definition 3.3.6}
Let the {\it $E$-dimension\/} of $G$ be the equivalence class
$\dim_E(G) = [S^{adG}] \in \Pic_E$ of the dualizing spectrum
$S^{adG}$ in the $E$-local Picard group. 
\enddefinition

\example{Example 3.3.7}
For $E = S$ the $S$-dimension of a compact Lie group $G$ equals its
manifold dimension in $\Pic_S \cong \Z$.  Similarly, for $E = H\F_p$
the $H\F_p$-dimension of a $p$-compact group $G$ is the same as its
$p$-dimension.
\endexample

\head 4. Computations \endhead

\subhead 4.1. A spectral sequence for $E$-homology \endsubhead

Suppose that the $S$-module $E$ is an $S$-algebra.
The standard left $G$-action $\alpha'$ on $DG_+$ makes $E_*(DG_+) =
E^{-*}(G)$ a left $E_*(G)$-module via the composite action map
$$
E_*(G) \otimes E_*(DG_+) \to E_*(S[G] \wedge DG_+)
@>\alpha'_*>> E_*(DG_+) \,.
$$

\proclaim{Proposition 4.1.1}
Let $E$ be an $S$-algebra and let $G$ be a stably dualizable group.
There is a spectral sequence
$$
E^2_{s,t} = \Tor^{E_*(G)}_{s,t}(E_*, E^{-*}(G))
\Longrightarrow
E_{s+t}(S^{-adG})
$$
converging strongly to $E_*(S^{-adG}) \cong E^{-*}(S^{adG})$.
\endproclaim

\demo{Proof}
This is the $E$-homology homotopy orbit spectral sequence, which is a
special case of the Eilenberg--Moore type spectral sequence \cite{EKMM97,
IV.6.4} for the $E$-homology of
$$
S^{-adG} = EG_+ \wedge_G DG_+ \cong S[EG] \wedge_{S[G]} DG_+
\,.
$$
(Other names in use are the bar spectral sequence and the
Steenrod--Rothenberg spectral sequence.)  Here $E_*(S[EG]) \cong E_*$,
$E_*(S[G]) \cong E_*(G)$ and $E_*(DG_+) \cong E^{-*}(G)$.  The duality
$S^{-adG} \simeq DS^{adG}$ from Proposition~3.2.3 relates the abutment
to the $E$-cohomology of $S^{adG}$.
\qed
\enddemo

\subhead 4.2. Morava $K$-theories \endsubhead

In this and the following section (4.3) we specialize to the case when
$E = K(n)$, for some fixed prime $p$ and number $0 \le n \le \infty$.
Hence stably dualizable means $K(n)$-locally stably dualizable, etc.

\proclaim{Lemma 4.2.1}
Let $G$ be a stably dualizable group, so that $H =
K(n)_*(G)$ is a finitely generated (free) module over $R = K(n)_*$.  Then
$H$ is a graded cocommutative Hopf algebra over $R$, and its $R$-dual $H^*
= K(n)^*(G) = \Hom_R(H, R)$ is a graded commutative Hopf algebra over $R$.
\endproclaim

\demo{Proof}
By \cite{HSt99, 8.6}, a topological group $G$ is stably dualizable
if and only if $H = K(n)_*(G)$ is finitely generated over $R = K(n)_*$.
The group multiplication and diagonal map on $G$ induce the Hopf algebra
structure on $H$, in view of the K{\"u}nneth isomorphism
$$
K(n)_*(X) \otimes_{K(n)_*} K(n)_*(Y) @>\cong>>
K(n)_*(X \wedge Y)
$$
in the case $X = Y = S[G]$.  The identity $K(n)^*(G) \cong
\Hom_R(H, R)$ is a case of the universal coefficient theorem
$$
K(n)^*(X) @>\cong>> \Hom_{K(n)_*}(K(n)_*(X), K(n)_*) \,.
$$
This also leads to the Hopf algebra structure on $H^*$.
\qed
\enddemo

\proclaim{Proposition 4.2.2}
Let $G$ be a stably dualizable group.  Then $K(n)_*(S^{adG}) \cong
\Sigma^d R$ for some integer $d$, and $K(n)_*(S^{-adG}) \cong \Sigma^{-d}
R$.  \endproclaim

\demo{Proof}
By Theorem~3.3.4, $S^{adG}$ is an invertible $K(n)$-local spectrum with
inverse $S^{-adG}$, so by the K{\"u}nneth theorem
$$
K(n)_*(S^{adG}) \otimes_R K(n)_*(S^{-adG}) \cong K(n)_*(S) = R
\,.
$$
This implies that $K(n)_*(S^{adG})$ and $K(n)_*(S^{-adG})$ both have
rank one over $R$.  (Alternatively, use Theorem~3.1.4 and the
K{\"u}nneth theorem to obtain the isomorphism
$$
H^* \otimes_R K(n)_*(S^{adG}) \cong H \,.
$$
The total ranks of $H^*$ and $H$ as $R$-modules are equal, and finite,
so $K(n)_*(S^{adG})$ must have rank one.  In view of \cite{HMS94, 1.3} or
\cite{HSt99, 14.2}, this also provides an alternative proof that $S^{adG}$
is invertible in the $K(n)$-local category.)
\qed
\enddemo

\definition{Definition 4.2.3}
Let the integer $d = \deg_{K(n)}(G)$ such that $K(n)_*(S^{adG}) \cong
\Sigma^d R$ be the {\it $K(n)$-degree\/} of $G$.  When $0 <
n < \infty$ this number is only well-defined modulo $|v_n| = 2(p^n-1)$.
\enddefinition

\remark{Remark 4.2.4}
The evident homomorphism $\deg \: \Pic_{K(n)} \to \Z/|v_n|$ takes the
$K(n)$-dimension of $G$ to its $K(n)$-degree.  By \cite{HMS94, 1.3} or
\cite{HSt99, 14.2} we also have $\widehat{E(n)}^*(S^{adG}) \cong
\Sigma^d \widehat{E(n)}^*$, where $\widehat{E(n)} = L_{K(n)} E(n)$.
Similarly $E_n^*(S^{adG}) \cong \Sigma^d E_n^*$, where $E_n$ is the
Hopkins--Miller commutative $S$-algebra.  Taking into account the
action of the $n$-th Morava stabilizer group on $E_n^*(S^{adG})$ it is
in principle possible to recover much more information about the
$K(n)$-dimension of $G$ than just the $K(n)$-degree.
\endremark

\medskip

For any graded commutative ring $R$ and $R$-algebra $H$, we may consider
both $H$ and its $R$-dual $H^* = \Hom_R(H, R)$ as left $H$-modules in
the standard way.  Recall from e.g.~\cite{Pa71, \S4} that $H$ is called
a {\it (graded) Frobenius algebra\/} over $R$ if
\roster
\item
$H$ is finitely generated and projective as an $R$-module, and
\item
$H$ and some suspension $\Sigma^d H^*$ are isomorphic as left $H$-modules.
\endroster
It follows that $H$ is also isomorphic to $\Sigma^d H^*$ as right
$H$-modules, and conversely.  A (left or right) module $M$ over a
Frobenius algebra $H$ is projective if and only if it is injective.

\proclaim{Proposition 4.2.4}
Let $G$ be a stably dualizable group.  Then $H = K(n)_*(G)$ is a
Frobenius algebra over $R = K(n)_*$.  In particular, $H^* = K(n)^*(G)$
is an injective and projective (left) $H$-module.  In fact, it is free
of rank one.
\endproclaim

\demo{Proof}
Applying $K(n)$-homology to the equivalence of Theorem~3.1.4 gives
an isomorphism
$$
H^* \otimes_R \Sigma^d R = K(n)_*(DG_+) \otimes_{K(n)_*} K(n)_*(S^{adG})
\cong K(n)_*(S[G]) = H \,.
$$
Here $H$ acts from the left via the inverse of the second $G$-action,
i.e., by the standard left action on $H^*$, the trivial action on
$K(n)_*(S^{adG}) = \Sigma^d R$, and the left action through inverses
on $H$.  We continue with the isomorphism
$$
\chi_* \: H = K(n)_*(G) @>\cong>> K(n)_*(G) = H
$$
induced by the conjugation $\chi$ on $S[G]$, which takes the left
action through inverses to the standard left action.  Then the
composite of these two isomorphisms exhibits $H$ as a Frobenius
algebra over $R$.

It is a formality that $H^*$ is injective as a left $H$-module, so
the general theory implies that it is also projective.  But we can also
see this directly in our case, since $H^* \cong \Sigma^{-d} H$
is an isomorphism of left $H$-modules, and the right hand side
is free of rank one and thus obviously projective.
\qed
\enddemo

\proclaim{Theorem 4.2.5}
Let $G$ be a $K(n)$-locally stably dualizable group.  The spectral
sequence
$$
E^2_{s,t} = \Tor^H_{s,t}(R, H^*)
\Longrightarrow
K(n)_{s+t}(S^{-adG})
$$
collapses to the line $s=0$ at the $E^2$-term.  The natural map $i \:
DG_+ \to S^{-adG}$ identifies
$$
\Sigma^{-d} R = K(n)_*(S^{-adG}) \cong R \otimes_H H^*
$$
with the left $H = K(n)_*(G)$-module indecomposables of $H^* =
K(n)_*(DG_+) = K(n)^{-*}(G)$.  Dually, the natural map $p \: S^{adG}
\to S[G]$ identifies
$$
\Sigma^d R = K(n)_*(S^{adG}) \cong \Hom_H(H^*, R)
$$
with the left $H^*$-comodule primitives in $H$.
\endproclaim

\demo{Proof}
The spectral sequence is that of Proposition~4.1.1 in the special case
$E = K(n)$.  By Proposition~4.2.4, $H^*$ is a free left $H$-module of
rank one, hence flat.  Thus $\Tor^H_{s,t}(R, H^*) = 0$ for $s > 0$,
while for $s=0$, $\Tor^H_{0,*}(R, H^*) = R \otimes_H H^*$.  Hence the
spectral sequence collapses to the line $s=0$, and the edge homomorphism
corresponding to the inclusion $i \: DG_+ \to EG_+ \wedge_G DG_+ =
S^{-adG}$ is the surjection $H^* = K(n)_*(DG_+) \to K(n)_*(S^{-adG})
= R \otimes_H H^*$.  Thinking of $H^*$ as a left $H$-module, these are
the $H$-module coinvariants, or indecomposables, of $H^*$.

Passing to duals, the projection $p \: S^{adG} = F(EG_+, S[G]) \to
S[G]$ is functionally dual to the inclusion above, hence induces the
$R$-dual injection $\Hom_R(R \otimes_H H^*, R) \to \Hom_R(H^*, R)$ in
$K(n)$-homology.  Thus $K(n)_*(S^{adG})$ is identified with $\Hom_R(R
\otimes_H H^*, R) \cong \Hom_H(H^*, R)$, sitting inside $\Hom_R(H^*,
R) \cong H$.  The left $H$-module structure on $H^*$ dualizes to a
left $H^*$-comodule structure on $H$.  The inclusion $\Hom_H(H^*, R)
\to \Hom_R(H^*, R) \cong H$ then identifies $\Hom_H(H^*, R)$ with the
$H^*$-comodule primitives in $H$.
\qed
\enddemo

\remark{Remark 4.2.6}
We sometimes write $Q_H(H^*) = R \otimes_H H^*$ for the left $H$-module
indecomposables of $H^*$, and dually $P_{H^*}(H) = \Hom_H(H^*, R)$
for the left $H^*$-comodule primitives in $H$.  Then $K(n)_*(S^{-adG})
\cong Q_H(H^*)$ and $K(n)_*(S^{adG}) \cong P_{H^*}(H)$.

To be explicit, an element $x \in H \cong \Hom_R(H^*, R)$ lies
in $\Hom_H(H^*, R)$ if and only if $(y * \xi)(x) = \xi(xy)$ equals
$\epsilon(y) \xi(x) = \xi(x \epsilon(y))$ for each $y \in H$ and $\xi
\in H^*$.  Here $\epsilon \: H \to R$ is the augmentation.  This condition
is equivalent to asking that $xy = 0$ for each $y \in \ker(\epsilon)$,
i.e., $x \in H$ multiplies to zero with each element in the augmentation
ideal of $H$.  So $P_{H^*}(H)$ is the left annihilator ideal of
the augmentation ideal of $H$.
\endremark

\subhead 4.3. Eilenberg--Mac\,Lane spaces \endsubhead

We can make the identifications in Theorem~4.2.5 explicit in the cases
when $G = K(\Z/p, q)$ is an Eilenberg--Mac\,Lane space.  For $p$ an
odd prime the $K(n)$-homology $H = K(n)_* K(\Z/p, q)$ was computed by
Ravenel and Wilson in \cite{RaW80, 9.2}, as we now recall.

Writing $K(n)_* K(\Z, 2) \cong K(n)_* \{\beta_m \mid m\ge0\}$ with
$|\beta_m| = 2m$ there are classes $a_m \in K(n)_* K(\Z/p, 1)$ in degree
$|a_m| = 2m$ for $0 \le m < p^n$ such that the Bockstein map $K(\Z/p, 1)
\to K(\Z, 2)$ takes each $a_m$ to $\beta_m$.  Let $a_{(i)} = a_{p^i}$
in degree $|a_{(i)}| = 2 p^i$ for $0 \le i < n$.  The $q$-fold cup
product $K(\Z/p, 1) \wedge \dots \wedge K(\Z/p, 1) \to K(\Z/p, q)$ takes
$a_{(i_1)} \otimes \dots \otimes a_{(i_q)}$ to a class $a_I \in K(n)_*
K(\Z/p, q)$, where $I = (i_1, \dots, i_q)$ and $|a_I| = 2(p^{i_1} +
\dots + p^{i_q})$.

For $q=0$, $G = K(\Z/p, 0) = \Z/p$ is a finite group and not very special
to the $K(n)$-local situation.  For each $q > n$, $K(\Z/p, q)$ has the
$K(n)$-homology of a point.  The intermediate cases $0 < q \le n$
are more interesting.

For $0 < q < n$ there is an algebra isomorphism
$$
K(n)_* K(\Z/p, q) \cong \bigotimes_I K(n)_*[a_I]/(a_I^{p^{\rho(I)}})
\,,
$$
where $I = (i_1, \dots, i_q)$ ranges over all integer sequences with $0 <
i_1 < \dots < i_q < n$, and $\rho(I) = s + 1$ where $s \in \{0, 1, \dots,
q\}$ is maximal such that the final $s$-term subsequence
has the form
$$
(i_{q-s+1}, \dots, i_q) = (n-s, \dots, n-1) \,.
$$
Equivalently, $s$ is minimal such that $i_{q-s} < n-s-1$.

For $q=n$ there is an algebra isomorphism
$$
K(n)_* K(\Z/p, n) \cong K(n)_*[a_I]/(a_I^p + (-1)^n v_n a_I)
\,,
$$
where $I = (0, 1, \dots, n-1)$.  Here $|a_I| = 2(1 + p + \dots + p^{n-1})
= 2(p^n-1)/(p-1)$.

\proclaim{Proposition 4.3.1}
For $G = K(\Z/p, q)$ with $0 < q < n$, $K(n)_*(S^{adG})$ is generated
over $K(n)_*$ by the product $\pi = \prod_I a_I^{p^{\rho(I)}-1}$.
Its $K(n)$-degree is $0$ modulo~$2(p^n-1)$.
\endproclaim

\demo{Proof}
By Theorem~4.2.5 we identify $K(n)_*(S^{adG})$ with the left
$H^*$-comodule primitives in $H$, which by Remark~4.2.6 consists of the
elements of $H$ that multiply to zero with every element in the
augmentation ideal of $H$.  These are generated by the product $\pi$
above.  Its degree $\deg_{K(n)}(G) \equiv |\pi|$ can be computed by
grouping together the integer sequences with the same value of $\rho(I)
= s+1$:
$$
\align
|\pi| &= \sum_I 2 (p^{i_1} + \dots p^{i_q}) (p^{\rho(I)}-1) \\
&= \sum_{0 \le s \le q \atop 1 \le i_1 < \dots < i_{q-s} \le n-s-2}
2 (p^{i_1} + \dots + p^{i_{q-s}} + p^{n-s} + \dots + p^{n-1}) (p^{s+1}-1) \\
\allowdisplaybreak
&\equiv
\sum_{0 \le s \le q \atop s+2 \le j_{s+1} < \dots < j_q \le n-1}
2 (p^{j_{s+1}} + \dots + p^{j_q} + p^1 + \dots + p^s) \\
& \qquad - \sum_{0 \le s \le q \atop 1 \le i_1 < \dots < i_{q-s} \le n-s-2}
2 (p^{i_1} + \dots + p^{i_{q-s}} + p^{n-s} + \dots + p^{n-1}) \\
&=
\sum_{1 \le j_1 < \dots < j_q \le n-1} 2 (p^{j_1} + \dots + p^{j_q})
- \sum_{1 \le i_1 < \dots < i_q \le n-1} 2 (p^{i_1} + \dots + p^{i_q})
= 0
\endalign
$$
modulo $2(p^n-1)$.
\qed
\enddemo

\proclaim{Proposition 4.3.2}
For $G = K(\Z/p, n)$, $K(n)_*(S^{adG})$ is generated over $K(n)_*$ by $\pi
= a_I^{p-1} + (-1)^n v_n$.  Its $K(n)$-degree is $0$ modulo~$2(p^n-1)$.
\endproclaim

\demo{Proof}
In this case the primitives in $K(n)_*(S^{adG})$ are generated by
$$
\pi = a_I^{p-1} + (-1)^n v_n
$$
in degree $|v_n| = 2(p^n-1)$.  So also in this case $\deg_{K(n)}(G)
\equiv 0$.
\qed
\enddemo

\remark{Remark 4.3.3}
It would be interesting to produce non-integer elements in the
$K(n)$-local Picard group $\Pic_{K(n)}$ as the class $[S^{adG}]$ of
the dualizing spectrum of a $K(n)$-locally stably dualizable group~$G$.
The Eilenberg--Mac\,Lane examples above do not decisively produce any
such non-integer elements.  Together with Lemmas~2.3.2 and~3.1.6, and
\cite{HRW98}, this indicates that the required stably dualizable group
$G$ should not even be homotopy commutative.  This adds interest to the
construction suggested in Remark~2.4.8.
\endremark

\head 5. Norm and transfer maps \endhead

\subhead 5.1. Thom spectra \endsubhead

The Thom space of a $G$-representation $V$ is the reduced Borel
construction, or homotopy orbit space, $BG^V = EG_+ \wedge_G S^V$,
where $S^V$ is the representation sphere of $V$.  Generalizing the
compact Lie group case, when $S^{adG}$ is the (suspension spectrum of)
the representation sphere of the adjoint representation $ad G$, we make
the following definition:

\definition{Definition 5.1.1}
Let $G$ be a stably dualizable group.  The {\it Thom spectrum}
$BG^{adG}$ of its dualizing spectrum is the homotopy orbit spectrum
$$
BG^{adG} = (S^{adG})_{hG} = EG_+ \wedge_G S^{adG} \,.
$$
The inclusion $G \subset EG$ induces the {\it bottom cell
inclusion} $i \: S^{adG} \to BG^{adG}$.
\enddefinition

Note that for $G$ abelian, $S^{adG}$ is a spectrum with
$EG$-action by Lemma~2.5.4, so in these cases
$$
BG^{adG} \simeq BG_+ \wedge S^{adG} \,.
$$

As in Proposition~4.1.1, when $E$ is an $S$-algebra there is a strongly
convergent spectral sequence
$$
\align
E^2_{s,t} &= \Tor^{E_*(G)}_{s,t}(E_*, E_*(S^{adG})) \\
&\Longrightarrow E_{s+t}(BG^{ad G}) \,.
\tag 5.1.2
\endalign
$$
When $E = K(n)$ we have $K(n)_*(S^{adG}) \cong \Sigma^d K(n)_*$
by Proposition~4.2.2, with $d = \deg_{K(n)}(G)$.  Thus the
spectral sequence takes the form
$$
\align
E^2_{s,t} &= \Tor^{K(n)_*(G)}_{s,t}(K(n)_*, \Sigma^d K(n)_*) \\
&\Longrightarrow K(n)_{s+t}(BG^{ad G}) \,.
\tag 5.1.3
\endalign
$$
When $S^{adG}$ is $K(n)$-orientable, so that the bottom cell inclusion
$i \: S^{adG} \to BG^{adG}$ induces a nonzero homomorphism $i_* \:
\Sigma^d K(n) \cong K(n)_*(S^{adG}) \to K(n)_*(BG^{adG})$, then this
spectral sequence~5.1.3 is a free comodule over the corresponding bar
spectral sequence for $K(n)_*(BG)$, on a single generator in degree $d$.

\subhead 5.2. The norm map and Tate cohomology \endsubhead

Let $X$ be a spectrum with left $G$-action.  We give it the trivial right
$G$-action.  The smash product $X \wedge S[G]$ then has the diagonal left
$G$-action, as well as the right $G$-action that only affects $S[G]$.
Consider forming homotopy orbits with respect to the left action,
and forming homotopy fixed points with respect to the right action.

We shall construct the norm map in three steps.
First, there is a canonical colimit/limit exchange map
$$
\kappa \: ((X \wedge S[G])^{hG})_{hG} @>>> ((X \wedge S[G])_{hG})^{hG}
\tag 5.2.1
$$
induced by the familiar map
$$
EG_+ \wedge F(EG_+, Y) @>>> F(EG_+, EG_+ \wedge Y) \,,
$$
in the case $Y = X \wedge S[G]$.

Second, there is a natural map $\nu \: X \wedge S^{adG} = X \wedge
S[G]^{hG} \to (X \wedge S[G])^{hG}$, since $G$ acts trivially on $X$
from the right.  It can be identified with the chain of weak equivalences
$$
X \wedge S^{adG} @>\simeq>> F(G_+, X \wedge S^{adG})^{hG}
@<\nu^{hG}<\simeq< (X \wedge DG_+ \wedge S^{adG})^{hG}
@>\simeq>> (X \wedge S[G])^{hG}
$$
where the middle map uses that $G$ is stably dualizable and the right
hand map uses the Poincar{\'e} duality equivalence of Theorem~3.1.4.
In particular $\nu \: X \wedge S^{adG} \to (X \wedge S[G])^{hG}$ is a
weak equivalence, and it induces a weak equivalence
$$
\nu_{hG} \: (X \wedge S^{adG})_{hG} @>\simeq>> ((X \wedge S[G])^{hG})_{hG}
$$
on homotopy orbits with respect to the left actions.  Note that $\nu_{hG}$
maps to the left hand side of~(5.2.1).  In the special case $X = S$,
the maps $\nu$ and $\nu_{hG}$ are isomorphisms.

Third, there is an untwisting equivalence $\zeta \: (X
\wedge S[G])_{hG} \to X \wedge S[G]_{hG} \simeq X \wedge S \cong X$,
cf.~\cite{LMS86, p.~76}, that takes the remaining right action on $(X
\wedge S[G])_{hG}$ to the right action on $X$ through the inverse of
the left action.  Hence there is an equivalence
$$
\zeta^{hG} \: ((X \wedge S[G])_{hG})^{hG} @>\simeq>> X^{hG}
$$
of homotopy fixed points, formed with respect to these right actions.
Note that $\zeta^{hG}$ maps from the right hand side of~(5.2.1).

\definition{Definition 5.2.2}
Let $X$ be a spectrum with left $G$-action.
The {\it (homotopy) norm map}
$$
N \: (X \wedge S^{adG})_{hG} @>>> X^{hG}
$$
is the composite of the natural maps:
$$
(X \wedge S^{adG})_{hG} @>\nu_{hG}>\simeq>
((X \wedge S[G])^{hG})_{hG} @>\kappa>> ((X \wedge S[G])_{hG})^{hG}
@>\zeta^{hG}>\simeq> X^{hG} \,,
$$
where $\nu_{hG}$ and $\zeta^{hG}$ are weak equivalences.
The {\it $G$-Tate cohomology spectrum} $X^{tG}$ of $X$ is the cofiber of the
norm map:
$$
(X \wedge S^{adG})_{hG} @>N>> X^{hG} \to X^{tG} \,.
$$
\enddefinition

\remark{Remark 5.2.3}
(a)
Note that $X^{tG} \simeq *$ if and only if the norm map $N$ is a weak
equivalence, which in turn is equivalent to the canonical colimit/limit
exchange map $\kappa$ in~(5.2.1) being a weak equivalence.  So $G$-Tate
cohomology measures the failure of $G$-homotopy orbits and $G$-homotopy
fixed points to commute.

(b)
In view of \cite{GM95, 3.5}, it is reasonable to expect that if $X$
is an $S$-algebra with $G$-action, then $X^{tG}$ is an $S$-algebra
and $X^{hG} \to X^{tG}$ is a map of $S$-algebras.  We do not know how
to give a direct model for $X^{tG}$, say as the ``$G$-fixed points''
of the spectrum $\widetilde{EG} \wedge F(EG_+, X)$ with $G$-action,
so it is not so easy to verify our expectation.  Here, as usual in this
context, $\widetilde{EG}$ is the mapping cone of the collapse map $c \:
EG_+ \to S^0$.
\endremark

\medskip

In the special case when $X = S$ with trivial $G$-action, the norm map
simplifies to the canonical colimit/limit exchange map
$$
BG^{adG} = (S^{adG})_{hG} = (S[G]^{hG})_{hG} @>\kappa>> (S[G]_{hG})^{hG}
\simeq S^{hG} \cong D(BG_+) \,.
$$
Here we use that $S^{hG} = F(EG_+, S)^G \cong F(BG_+, S) = D(BG_+)$
is the functional dual of $BG_+$, since $S$ has trivial $G$-action.
Hence there is a cofiber sequence
$$
BG^{adG} @>N>> D(BG_+) @>>> S^{tG} \,.
$$

In the case of a compact Lie group $G$, the $G$-Tate cohomology $X^{tG}$
is the same as that denoted $t_G(X)^G$ by Greenlees and May \cite{GM95}
and $\hat\Bbb H(G, X)$ by B{\"o}kstedt and Madsen \cite{BM94}.

\definition{Definition 5.2.4}
A spectrum with $G$-action $X$ is in the {\it thick subcategory\/}
generated by spectra of the form $G_+ \wedge W$, if $X$ can be built from
$*$ in finitely many steps by (1) attaching cones on induced $G$-spectra
of the form $G_+ \wedge W$, with $W$ any spectrum, (2) passage to
(weakly) equivalent spectra with $G$-action and (3) passage to retracts.
For instance, any finite $G$-cell spectrum has this form.
\enddefinition

\proclaim{Theorem 5.2.5}
Let $G$ be a stably dualizable group.  If a spectrum with $G$-action
$X$ is in the thick subcategory generated by spectra of the form $G_+
\wedge W$, then:
\roster
\item
The norm map $N \: (X \wedge S^{adG})_{hG} \to X^{hG}$ for $X$ is an
equivalence.
\item
The $G$-Tate cohomology $X^{tG} \simeq *$ is contractible.
\endroster
\endproclaim

\demo{Proof}
If $X = G_+ \wedge W$ is induced up from a spectrum $W$ with trivial
$G$-action, the source of the norm map can be identified with
$$
(G_+ \wedge W \wedge S^{adG})_{hG} \simeq W \wedge S^{adG}
$$
while the target of the norm map can be identified with
$$
(G_+ \wedge W)^{hG} \simeq (DG_+ \wedge S^{adG} \wedge W)^{hG} \simeq
F(G_+, W \wedge S^{adG})^{hG} \simeq W \wedge S^{adG}
\,.
$$
These identifications are compatible, as can be checked by starting
with the case $W = S$, hence in this case the norm map is itself an
equivalence.  The general case follows by induction on the number of
attachments made.
\qed
\enddemo

\remark{Remark 5.2.6}
This result generalizes the third case of \cite{Kl01, Thm.~D}, from
compact Lie groups to stably dualizable groups.  For compact Lie
groups $G$ this norm equivalence can be compared with the genuinely
$G$-equivariant Adams equivalence $Y/G \simeq (Y \wedge S^{-adG})^G$
for $Y$ a free $G$-spectrum \cite{LMS86, II.7}.  Any such $Y$ is a filtered
colimit of finite, free $G$-spectra, which are in the thick subcategory
generated by $S[G] \cong G_+ \wedge S$.  But, while genuine $G$-fixed points
($Y \mapsto Y^G$) commute with filtered colimits, this is not generally
the case for $G$-homotopy fixed points ($Y \mapsto Y^{hG}$).  Therefore we
cannot extend Theorem~5.2.5 to all spectra $X$ with free $G$-action.
\endremark

\medskip

There is also a dual construction $X_{tG}$ that is to Tate homology as
the Tate construction $X^{tG}$ is to Tate cohomology.  To define it,
we suppose that $X$ is a spectrum with right $G$-action, and give it the
trivial left $G$-action.  The smash product $DG_+ \wedge X$ has the left
$G$-action that only affects $DG_+$, and the diagonal right $G$-action.
There is a canonical colimit/limit exchange map
$$
\kappa \: ((DG_+ \wedge X)^{hG})_{hG} @>>> ((DG_+ \wedge X)_{hG})^{hG}
\,.
$$
The source of $\kappa$ receives an equivalence from $X_{hG}$ obtained
by applying homotopy orbits to the equivalence $X \to F(G_+, X)^{hG}$.
The target of $\kappa$ admits a weak equivalence to $(S^{-adG} \wedge
X)^{hG}$ obtained by taking homotopy fixed points of the isomorphism
$$
(DG_+ \wedge X)_{hG} @>\cong>> (DG_+)_{hG} \wedge X = S^{-adG} \wedge X
\,,
$$
which exists because $G$ acts trivially on $X$ from the left.

\definition{Definition 5.2.7}
Taken together, these maps yield the alternate norm map $N'$, defined
as the composite:
$$
N' \: X_{hG} @>\simeq>> ((DG_+ \wedge X)^{hG})_{hG} @>\kappa>>
((DG_+ \wedge X)_{hG})^{hG} @>\simeq>> (S^{-adG} \wedge X)^{hG} \,.
$$
Its homotopy fiber $X_{tG}$ is the {\it $G$-Tate homology spectrum},
and sits in a cofiber sequence
$$
X_{tG} @>>> X_{hG} @>N'>> (S^{-adG} \wedge X)^{hG} \,.
$$
\enddefinition

If $X$ is an $S$-coalgebra with $G$-action, it again appears likely
that $X_{tG}$ is such a coalgebra and that $X_{tG} \to X_{hG}$ is a map
of $S$-coalgebras.

\subhead 5.3. The $G$-transfer map \endsubhead

\definition{Definition 5.3.1}
Let $X$ be a spectrum with left $G$-action.
The {\it $G$-transfer map}
$$
\trf_G \: (X \wedge S^{adG})_{hG} @>>> X
$$
is the composite of the norm map $N \: (X \wedge S^{adG})_{hG} \to X^{hG}$
and the forgetful map $p \: X^{hG} \to X$.

When $X = S[Y]$ is the unreduced suspension spectrum of a $G$-space $Y$,
this is the {\it dimension-shifting $G$-transfer map\/} associated to
the principal $G$-bundle $Y \simeq EG \times Y \to EG \times_G Y$.
\enddefinition

\subhead 5.4. $E$-local homotopy classes \endsubhead

Let $G$ be a stably dualizable group, with dualizing spectrum $S^{adG}$.
The composite of the bottom cell inclusion $i \: S^{adG} \to BG^{adG}$
and the dimension-shifting $G$-transfer $\trf_G \: BG^{adG} \to S$
is the composite map
$$
S^{adG} @>i>> (S^{adG})_{hG} = (S[G]^{hG})_{hG}
@>\kappa>> (S[G]_{hG})^{hG} \simeq S^{hG} @>p>> S \,.
$$
Noting that the projection $p$ amounts to forgetting $G$-homotopy
invariance, this map can also be expressed as the composite
$$
S^{adG} = S[G]^{hG} @>p>> S[G] @>i>> S[G]_{hG} \simeq S
\,.
$$

\definition{Definition 5.4.1}
The composite map $p \kappa i \: S^{adG} \to S$ represents a class
denoted $[G] \in \pi_*(L_E S)$ in the $\Pic_E$-graded homotopy groups
of the $E$-local sphere spectrum, in degree $* = \dim_E(G) = [S^{adG}]
\in \Pic_E$.  We might call $[G]$ the $E$-local {\it stably framed
bordism class\/} of $G$.
\enddefinition

\example{Example 5.4.2}
For the circle group $G = S^1$ and $E=S$ we have $[G] = \eta \in
\pi_1(S)$.  For the $p$-complete Sullivan sphere $G = (S^{2p-3})^\wedge_p$
and $E = H\F_p$ we have $[G] = \alpha_1 \in \pi_{2p-3}(S^\wedge_p)$,
when $p$ is odd.  These examples are also detected $K(1)$-locally, i.e.,
for $G$ considered as a $K(1)$-locally stably dualizable group.
\endexample

\proclaim{Lemma 5.4.3}
In the case $E = K(n)$, the induced homomorphism
$$
[G]_* \: \Sigma^d K(n)_* \cong K(n)_*(S^{adG}) \to K(n)_*(S) = K(n)_*
$$
takes a generator of the $K(n)^*(G)$-comodule primitives in
$K(n)_*(G)$ to its image under the augmentation
$\epsilon \: K(n)_*(G) \to K(n)_*$.
\endproclaim

\demo{Proof}
Recall from Theorem~4.2.5 that $K(n)_*(S^{adG})$ is identified with the
$H^* = K(n)^*(G)$-comodule primitives $\Hom_H(H^*, R)$, and the
projection $p \: S^{adG} \to S[G]$ induces the forgetful inclusion into
$\Hom_R(H^*, R) \cong H = K(n)_*(G)$.  The inclusion $i \: S[G] \to
S[G]_{hG} \simeq S$ induces the augmentation $\epsilon$, which
establishes the claim
\qed
\enddemo

\example{Example 5.4.4}
When $E = K(n)$ and $G$ is a finite discrete group we get $H = R[G]$
and $P_{H^*}(H) \cong R\{N\}$, where $N = \sum_{g \in G} g$ is the
norm element in $H$.  Then $\epsilon(N) = |G|$ equals the order of $G$,
so $[G]_*$ multiplies by $|G|$ in $R = K(n)_*$.
\endexample

\example{Example 5.4.5}
When $E = K(n)$ and $G = K(\Z/p, q)$ for $0 < q < n$, the
$H^*$-comodule primitives were found in Proposition~4.3.1 to be
generated by an element $\pi$ that lies in the augmentation ideal
$\ker(\epsilon)$, so the induced homomorphism $[G]_*$ is zero and $[G]
\: S^{adG} \to S$ has positive $K(n)$-based Adams filtration.
\endexample

\example{Example 5.4.6}
When $E = K(n)$ and $G = K(\Z/p, n)$, with $q=n$, Proposition~4.3.2
exhibited a generating element $\pi = a_I^{p-1} + (-1)^n v_n$ for
$P_{H^*}(H)$, which augments to the unit $(-1)^n v_n \in K(n)_*$.
Hence in this case $[G] \: S^{adG} \to S$ induces an isomorphism on
$K(n)$-homology, and so $S^{adG} \simeq S$ in the $K(n)$-local
category.  By Lemma~2.5.4, the $G$-action on $S^{adG}$ is homotopy
trivial in this case.  Hence the Poincar{\'e} duality equivalence~3.1.4
amounts to a $K(n)$-local self-duality equivalence
$$
F(G_+, L_{K(n)}S) = DG_+ \simeq S[G] = L_{K(n)} \Sigma^\infty G_+
\tag 5.4.7
$$
for $G = K(\Z/p, n)$, which is left and right $G$-equivariant up to
homotopy, and which may be compared with \cite{HSt99, 8.7}.
\endexample

\enddocument